\documentclass[review]{elsarticle}
\journal{Signal Processing}

\usepackage[utf8]{inputenc}
\usepackage[cmex10]{amsmath}
\usepackage{booktabs}
\usepackage{comment}

\newtheorem{lemma}{Lemma}

\newtheorem{definition}{Definition}

\newproof{IEEEproof}{Proof}

\usepackage[algoruled]{algorithm2e}
\SetKwInput{KwParam}{Parameters}
\SetKwInput{KwNotes}{Notes}
\SetKwComment{Comment}{}{}
\SetCommentSty{small}
\makeatletter
\newcommand{\removelatexerror}{\let\@latex@error\@gobble}
\makeatother

\usepackage{pgfplots}
\pgfplotsset{compat=1.14}
\usetikzlibrary{plotmarks}
\newlength\figureheight 
\newlength\figurewidth 

\usepackage[a4paper, margin=1.25in]{geometry}

\usepackage{xspace}		
\usepackage{amssymb}	
\usepackage{mathtools}	
\usepackage{bbold}		

\makeatletter
 \begingroup
 \catcode`\_=\active
 \protected\gdef_{\@ifnextchar|\subtextup\sb}
 \endgroup
 \def\subtextup|#1|{\sb{\textup{#1}}}
 \AtBeginDocument{\catcode`\_=12 \mathcode`\_=32768 }
\makeatother

\newcommand{\shownumber}{\addtocounter{equation}{1}\tag{\theequation}}
\newcommand{\ii}{^{-1}}
\newcommand{\db}{\ensuremath{\,\textrm{dB}}\xspace}

\newcommand{\intr}{\mathop{\bf int}}
\newcommand{\conv}{\mathop{\bf conv}}
\newcommand{\bdy}{\mathop{\bf bd}}

\newcommand{\h}{^\mathrm{H}}
\newcommand{\T}{^{\mathrm{T}}}

\DeclareMathOperator{\diag}{diag}

\DeclareMathOperator{\tr}{tr}

\DeclareMathOperator{\re}{Re}
\DeclareMathOperator{\im}{Im}

\newcommand{\E}[1]{\mathbb{E}\!\left[#1\right]}
\newcommand{\norm}[1]{\left\lVert#1\right\rVert}

\newcommand{\bbr}{\mathbb{R}}
\newcommand{\bbc}{\mathbb{C}}

\newcommand{\cala}{\ensuremath{\mathcal{A}}\xspace}

\newcommand{\calc}{\ensuremath{\mathcal{C}}\xspace}

\newcommand{\calk}{\ensuremath{\mathcal{K}}\xspace}

\newcommand{\calm}{\ensuremath{\mathcal{M}}\xspace}

\newcommand{\calo}{\ensuremath{\mathcal{O}}\xspace}

\usepackage[final,citecolor=black,bookmarks,bookmarksnumbered,colorlinks=true,linkcolor=black]{hyperref}

\begin{document}

\begin{frontmatter}

\title{A Fast Interior-Point Method for Atomic Norm Soft Thresholding}

\author{Thomas Lundgaard Hansen\fnref{tlh}}
\fntext[tlh]{The work of T.~L.~Hansen is supported by the Danish Council for
	Independent Research under grant id DFF--4005--00549. The work of T.~L.~Jensen was partly supported by by Danich Council for Independent Research under grant id DFF--4005--00122}

\author{Tobias Lindstrøm Jensen\corref{mycorrespondingauthor}}
\cortext[mycorrespondingauthor]{Corresponding author}
\ead{tlj@its.aau.dk}

\address{
	Aalborg University, Fr. Bajersvej 7, DK-9220 Aalborg, Denmark
}


\begin{abstract}
	The atomic norm provides a generalization of the $\ell_1$-norm
        to continuous parameter spaces. When applied as a sparse
        regularizer for line spectral estimation the solution can be
        obtained by solving a convex optimization problem. This
        problem is known as atomic norm soft thresholding (AST). It
        can be cast as a semidefinite program and solved by standard
        methods. In the semidefinite formulation there are
        $\calo(N^2)$ dual variables which complicates the implementation of a standard primal-dual
        interior-point method based on symmetric cones. That has lead researchers
        to consider the alternating direction method of multipliers
        (ADMM) for the solution of AST, but this method is still
        somewhat slow for large problem sizes.  To obtain a faster
        algorithm we reformulate AST as a non-symmetric conic
        program. That has two properties of key importance to its
        numerical solution: the conic formulation has only $\calo(N)$
        dual variables and the Toeplitz structure inherent to AST is
        preserved. Based on it we derive FastAST which is a
        primal-dual interior-point method for solving AST. Two
        variants are considered with the fastest one requiring only
        $\calo(N^2)$ flops per iteration. Extensive numerical
        experiments demonstrate that both variants of FastAST solve
        AST significantly faster than a state-of-the-art solver based
        on ADMM.
\end{abstract}

\begin{keyword}
	Atomic Norm Minimization \sep
	Atomic Norm Soft Thresholding \sep
	Line Spectral Estimation \sep
	Convex Optimization \sep
	Interior-Point Methods \sep
	Non-symmetric Conic Optimization
\end{keyword}

\end{frontmatter}


\section{Introduction}

It is well known that sparse estimation problems can be formulated as convex
optimization problems using the $\ell_1$-norm. The $\ell_1$-norm can be
generalized to continuous parameter spaces through the so-called atomic norm
\cite{Chandrasekaran:2012}.
Convex modelling of sparsity constraints has two highly attractive traits: 
convex optimization problems can easily be solved both in theory \cite{Ne:94} and in
practice \cite{Wri:97,Nocedal99}, and, a number of recovery
guarantees can be obtained within this framework. Such recovery guarantees are
studied in signal processing under the name compressed sensing \cite{Candes:06, Donaho:06,
candes:2008} and they generalize nicely to the atomic norm minimization
approach \cite{Tang:2013,Bhaskar:2013,Tang:2015,Candes:2014,Candes:2013}.

The most prominent example of estimation with the atomic norm is the
application to line spectral estimation
\cite{Bhaskar:2013,Candes:2014,Candes:2013}, in which case it is known as
atomic norm soft thresholding (AST).
The popularity of AST is, partly, due to the fact that it can be cast as a
semidefinite programming (SDP) problem (we refer to Sec.~\ref{sec:review} for
a review of  AST,)
\begin{equation}
\begin{array}{ll}
	\text{minimize}_{v,x,u} &
	\norm{x-y}_2^2 + \tau (v + w\T u) \\[1mm]
	\text{subject to} &
	\left(\begin{matrix}
		T(u) & x \\ x\h & v
	\end{matrix}\right)
	\succeq 0,
\end{array}
\label{problem}
\end{equation}
where $v\in\bbr$, $x\in\bbc^N$, $u\in\bbr^{2N-1}$ are the variables of the
problem and $y\in\bbc^N$, $\tau\in\bbr$, $w\in\bbr^{2N-1}$ are fixed (known)
parameters.
The function $T(u):\bbr^{2N-1}\rightarrow\bbc^{N\times N}$ outputs a
complex Hermitian Toeplitz matrix constructed from
$u$, such that the first row is $(2u_0,\ldots,u_{N-1}) +
j(0,u_N,\ldots,u_{2N-2})$. To be precise, AST is obtained by selecting
$w=2e_0$ in \eqref{problem}, where $e_0$ is a vector with $1$ in the first entry and zeros
elsewhere.

The state-of-the-art method for solving \eqref{problem} is via the
alternating direction method of multipliers (ADMM) and used in
\cite{Bhaskar:2013, Cho:2015, li-multiple, Yang:15}.
While this method is reasonably fast, it has some drawbacks. It requires the
calculation of an eigenvalue decomposition in each iteration at cost
$\calo(N^3)$ floating-point operations (flops). This means that for large $N$
it is exceedingly slow. As is often seen with proximal methods it also has
slow convergence if a solution of high accuracy is requested.

Da Costa \textit{et al.} \cite{Costa:2017} apply a low-dimensional projection of the
observation vector to reduce the problem size and therefore the computational
complexity of AST. In the noise-free case and under certain regularity
conditions, it is shown that the estimation accuracy is not affected by doing
so. However, it is clear that this approach discards observed data and the
estimation accuracy will be degraded in the noisy case. Another attempt at a
fast solver for AST is \cite{Rao:2015}, but due to its non-SDP implementation
utilizing a frequency grid and real, positive coefficients, this approach only
allows for a fixed known phase for all components and hence cannot solve the
line spectral problem as described in this paper. There has been other attempts
to solve atomic norm problems efficiently in e.g.~\cite{Boyd:2017,
Vinyes:2017}, but without covering the AST problem.

The formulation of AST in \eqref{problem} is casted as an SDP problem. SDP
problems have been subject to intensive research since the 1990's and their
solution using primal-dual interior-point methods (IPMs) is now understood
well \cite{Ne:94, Wri:97, BoVa:04, Sturm:99, Toh:1999}.
The Lagrangian dual of \eqref{problem} has $\calo(N^2)$ dual
variables due to the semidefinite matrix constraint. The direct application of
a standard primal-dual IPM thus requires $\calo(N^6)$ flops per iteration at
best (using direct methods for solving linear systems of equations), but can be reduced by eliminating the dual variables from the linear system (in general, the exact number flops per iteration will depend on the implementation). Compared
to this approach, proximal methods (such as ADMM) which require $\calo(N^3)$
flops per iteration are preferable, even if they converge much slower than
primal-dual IPMs. That explains why primal-dual IPMs have
not gained traction for the
solution of \eqref{problem}. In this work we reformulate the constraint in
\eqref{problem} as a non-symmetric conic constraint on the vector
$(v,x\T,u\T)\T$. This formulation immediately reduces the number of dual
variables to $\calo(N)$ and sets the scene for a reintroduction of primal-dual
IPMs as a very competitive class of algorithms for solution of AST.

Primal-dual IPMs for conic programming typically rely on a symmetry between the
primal and dual problems. The formulation of
such symmetric primal-dual IPMs requires the existence of a
self-scaled barrier function for the cone involved in the constraint
\cite{NeTo-selfscaled, NeTo-primaldual}. Güler
\cite{guler-barrier}
showed that such barrier functions exist only for the class of homogeneous and
self-dual cones. The cone in our formulation is not self-dual and so a symmetric
primal-dual IPM cannot be formulated.
Non-symmetric conic optimization has received some attention
\cite{nesterov-towards, skajaa-implementing, skajaa-homogeneous,
tuncel-generalization}. These methods generally rely on the availability of a
barrier function for the dual cone and possibly evaluation of its gradient and
Hessian. An easy-to-calculate dual barrier is not available for the cone associated to
the constraint of our formulation; only an oracle which can determine
membership in the dual cone is available (part of our contribution is to show
how such an oracle can be constructed.)

To derive a non-symmetric primal-dual IPM which does not rely on evaluating the
dual barrier or its derivatives, we formulate the augmented Karush-Kuhn-Tucker
conditions and devise a dedicated approach to solving these. This approach is
shown to converge to a primal-dual feasible point. A lower bound on the
objective function is calculated in those iterations where a dual feasible
point (as determined by the oracle) is available. From the lower bound a
duality gap can be evaluated, thus providing a method for dynamically updating
the barrier parameter. We show that the proposed method enjoys global
convergence.

Our focus is on obtaining an algorithm which has fast runtime \textit{in
practice}, i.e., it has both low per-iteration computational complexity and it
exhibits reasonably fast convergence. Theoretical statements regarding
for example convergence speed are left for future work.
At the core of obtaining a practically fast algorithm lies the already
mentioned conic formulation (which brings the number of dual variables down to
$\calo(N)$), along with techniques for fast evaluation of linear algebra in
each step of the algorithm. These evaluations are based on fast algorithms
\cite{ammar-generalized, ammar-numerical, levinson-wiener, durbin-fitting} for
inversion of Toeplitz matrices. Related techniques are employed in
\cite{HansenFR17, musicus-fastmlm, genin-optimization, alkire-autocorrelation}.

We dub the algorithm FastAST.
Both Newton's method and a quasi-Newton method are considered for evaluation of
the search direction in FastAST. When using Newton's method the algorithm
requires $\calo(N^3)$ flops per iteration, while the quasi-Newton variant only
requires $\calo(N^2)$ flops per iteration. The numerical experiments in Sec.
\ref{sec:numerical} show that the quasi-Newton variant is faster in practice.
Due to numerical inaccuracies in the calculation of the search direction the
quasi-Newton variant is not able to obtain a solution of very high accuracy.
Solving \eqref{problem} to high accuracy makes a difference with very large
signal-to-noise ratios and in these cases the variant of FastAST using Newton's
method should be used. Both the Newton's and quasi-Newton variants of FastAST
are significantly faster than the ADMM-based solvers for \eqref{problem}.

Along with the primal-dual IPM presented here, we have also experimented with a
primal-only version which is simpler to derive. The primal-only
approach does not provide a good way to select the barrier parameter (which we
denote $t$, see Sec. \ref{sec:barriers}). This in turn forces the primal-only
approach to use either
overly conservative short-step \cite{nemirovski-lecturenotes} updates of the
barrier parameter or it requires the barrier problem to be solved to high
accuracy for each fixed $t$. Both scenarios lead to a significant increase in
the number of iterations required by the primal-only algorithm, resulting in
significantly increased runtime.  On the contrary, the primal-dual version
presented in this paper allows for evaluation of a duality gap in each
iteration. The duality gap gives a natural way to select the barrier parameter
and also provides a very precise stopping criterion.

The paper is outlined as follows. In Sec.~\ref{sec:review} we begin with a
brief review of atomic norm minimization and its application to line spectral
estimation. Sec.~\ref{sec:non-symmetric} details the formulation of
\eqref{problem} as a non-symmetric conic optimization program along with the
theory of its solution. In Sec.~\ref{sec:method} we present our numerical
algorithm along with implementation details. The exploitation of Toeplitz
structure for fast evaluation of each step in the algorithm is discussed in
Sec.~\ref{sec:complexity}. Numerical experiments which validate the
practical efficacy of the proposed algorithm are presented in
Sec.~\ref{sec:numerical}.

\section{A Review of Atomic Norm Soft Thresholding}
\label{sec:review}

\subsection{Line Spectral Estimation}

Consider an observation vector $y\in\bbc^N$,
\begin{align}
	y = x + \zeta,
	\qquad
	x = \sum_{k=0}^{K-1} c_k a(\omega_k),
	\label{model}
\end{align}
where $\zeta\in\bbc^N$ is a noise vector and $x\in\bbc^N$ is a signal of
interest composed of $K$ sinusoids, each with angular frequency
$\omega_k\in[0,2\pi)$ and complex coefficient $c_k\in\bbc$. The steering vector
$a(\omega)$ has
entries $\left(a(\omega)\right)_n = \exp(jn\omega)$ for $n=0,\ldots, N-1$ and
$j=\sqrt{-1}$. In line spectral estimation
the task is to estimate the values $(K, c_0, \ldots, c_{K-1}, \omega_0, \ldots,
\omega_{K-1})$. The crux of line spectral estimation lies in obtaining
estimates of the
model order $K$ and the frequencies $\{\omega_k\}$. When these are available
the coefficients $\{c_k\}$ can easily be estimated using a least-squares
approach.
The problem is ubiquitous in signal processing; examples include direction
of arrival estimation using sensor arrays \cite{malioutov-sparse,
ottersten-analysis}, bearing and range estimation in synthetic aperture radar
\cite{carriere-high}, channel estimation in wireless communications
\cite{bajwa-sparsechannel} and simulation of atomic systems in molecular
dynamics \cite{andrade-application}.

\subsection{Modelling Sparsity With the Atomic Norm}
The atomic norm \cite{Chandrasekaran:2012, Tang:2013, Bhaskar:2013,
Tang:2015} provides a tool for describing notions of sparsity in a general
setting. It is defined in terms of the \textit{atomic set} $\cala$.  Each
member of $\cala$ is referred to as an \textit{atom}. The atoms are the basic
building block of the signal and define the dictionary in which the signal has
a sparse representation. The atomic norm induced by $\cala$ is defined as
\begin{align}
	\lVert x \rVert_\cala = \inf \{ \alpha>0 : x \in \alpha \conv \cala \},
\end{align}
where $\conv\cala$ is the convex hull of $\cala$.

For line spectral estimation the atomic set is selected as the set of complex
rotated Fourier atoms \cite{ Tang:2013, Bhaskar:2013, Tang:2015}
\begin{align}
	\cala = \{ a(\omega) \exp(j\phi) : \omega\in[0,2\pi), \phi\in[0,2\pi) \}
\end{align}
and the corresponding atomic norm can be described as
\begin{align}
	\lVert x \rVert_\cala = \inf_{K, \{c_k, \omega_k\} }
		\left\{ \sum_{k=0}^{K-1} |c_k|
		: x = \sum_{k=0}^{K-1} c_k a(\omega_k) \right\}.
	\label{atomic_norm}
\end{align}
It is clear that the atomic norm provides a generalization of the $\ell_1$-norm
to the continuous parameter space $\omega_k\in[0,2\pi)$. Through the use of a
dual polynomial characterization the atomic norm can be expressed as the
solution of an SDP,
\begin{equation}
\begin{array}{rll}
	\lVert x \rVert_\cala = \hspace{-\medskipamount}
	& \text{minimize}_{v,u} & \frac{1}{2} \left(v + \frac{1}{N}\tr T(u)\right) \\[1mm]
	& \text{subject to} &
		\left(\begin{matrix}
			T(u) & x \\ x\h & v
		\end{matrix}\right)
		\succeq 0.
\end{array}
\label{atomic_sdp}
\end{equation}

\subsection{Atomic Norm Soft Thresholding}
AST \cite{Bhaskar:2013} is inspired by the least absolute shrinkage and selection operator (LASSO)
\cite{Tibshirani:94} and solves
\begin{equation}
	\label{ast}
\begin{array}{rll}
	\text{minimize}_{x} &
		\lVert x-y \rVert_2^2 + 2\tau\lVert x \rVert_\cala,
\end{array}
\end{equation}
where $\tau>0$ is a regularization parameter to be chosen. It is
clear that AST is recovered in \eqref{problem} by selecting $w=2e_0$.

Once a solution $(v^\star, x^\star, u^\star)$ of \eqref{ast} has been found,
estimates of the model order $K$ and the frequencies $\{\omega_k\}$ can be
obtained by examining a certain dual polynomial constructed from $x^\star$.
This process determines the solution in \eqref{atomic_norm} for the recovered
signal $x^\star$. Under a, somewhat restrictive, assumption concerning
separation of the frequencies $\{\omega_k\}$ a number of theoretical statements
can be given regarding signal and frequency recovery using AST. We refer to
\cite{Tang:2013, Bhaskar:2013, Tang:2015, Candes:2014, Candes:2013} for
details.

We now consider the selection of the regularization parameter $\tau$. Clearly
the choice of $\tau$ crucially influences the estimation accuracy of AST. It is
this parameter which determines the trade-off between fidelity and sparsity
which is inherent in any estimator involving the model order $K$. With all else
being equal, selecting larger $\tau$ gives estimates with smaller values of
$K$.
Let $\lVert \cdot \rVert_\cala^*$ denote the dual norm of the atomic norm $\lVert
\cdot\rVert_\cala$. Then the theoretical analysis in \cite{Bhaskar:2013}
requires $\tau \ge \E{\lVert \zeta \rVert_\cala^*}$. For a white, zero-mean circularly symmetric
complex Gaussian noise vector $\zeta$ with entry-wise variance $\sigma^2$ such an
upper bound is given by \cite{Bhaskar:2013},
\begin{align}
	\tau = \sigma
		\frac{\log(N)+1}{\log(N)}
		\sqrt{
			N\log(N) + N\log(4\pi\log(N))
			},
		\label{tau}
\end{align}
where $\log()$ is the natural logarithm. This choice has been shown to perform
well in practice and we also use it in our simulation study.

\section{Non-symmetric Conic Optimization}
\label{sec:non-symmetric}


We now return to our main focus: That of numerically solving the conic program \eqref{problem}.
It can be written in the form
\begin{equation}
	\label{conicproblem}
\begin{array}{ll}
	\text{minimize} & f(\mu) \\[1mm]
	\text{subject to} & \mu\in\calk,
\end{array}
\end{equation}
where $f(\mu) = \norm{x-y}_2^2+\tau(v+w\T u)$ and $\calk$ is the cone defined by
\begin{align}
	\calk \triangleq \left\{
		\mu = \left(\begin{matrix}
			v \\ x \\ u
		\end{matrix}\right) :
			\left(\begin{matrix}
				T(u) & x \\ x\h & v
			\end{matrix}\right)
			\succeq 0
		\right\}.
\end{align}
As a precursor to deriving a primal-dual IPM, we explore
the properties of $\calk$ and its dual.
It is easy to show that $\calk$ is a proper cone (convex, closed, solid and
pointed; see \cite{BoVa:04}).
The dual cone $\calk^*$ of $\calk$ is defined as
\begin{align}
	\calk^* = \left\{
		\lambda
			: \left< \lambda, \mu \right> \ge 0 \;\; \forall \; \mu\in\calk
		\right\}.
		\label{Kstar_def}
\end{align}
Since $\calk$ is proper, so is $\calk^*$ \cite{BoVa:04}.

In this paper (primal) variables in the cone $\calk$ are denoted by
$\mu=(v,x\T,u\T)\T$. The (dual) variables in the cone $\calk^*$ are denoted by
$\lambda=(\rho,s\T,z\T)\T$, with $\rho\in\bbr$, $s\in\bbc^{N}$ and $z\in\bbr^{2N-1}$.
The inner product between them is defined as
$\left<\lambda,\mu\right> = \rho v + \re(s\h x) + z\T u$.

\subsection{Checking for Dual Cone Membership}
\label{dualmembership}
In our proposed method, we need to check for $\lambda\in\calk^*$.
In order to characterize the dual cone $\calk^*$, the cone of positive semidefinite
Hermitian Toeplitz matrices is needed:
\begin{align}
	\calc \triangleq \{ u \in\bbr^{2N-1} : T(u) \succeq 0 \}.
\end{align}
This cone is also proper.
The corresponding dual cone $\calc^*$ is defined analogously to
\eqref{Kstar_def}.

Let the function $T^*$ be the adjoint%
\footnote{
	$T^*$ is easy to calculate:
	Let $B$ be Hermitian and let $\beta_n$ denote the sum over the
	$n$th upper diagonal of $B$,
	\begin{align*}
		\beta_n = \sum_{m=0}^{N-1-n} B_{m,m+n},
	\end{align*}
	for $n=0,\ldots,N-1$.
	Then\\
	$T^*\!(B)=(2\beta_0, 2\re(\beta_1), \dots, 2\re(\beta_{N-1}),
	2\im(\beta_1), \dots, 2\im(\beta_{N-1}) )\T$.
}
of the linear map $T$, i.e., $T^* : \bbc^{N\times N} \rightarrow \bbr^{2N-1}$
is such that for every Hermitian $B\in\bbc^{N\times N}$ we have $\tr(T(u)\h B)
= T^*(B)\T u$.
Then we then have the following lemma.
\begin{lemma}
	\label{Kstar_lemma}
	The dual cone of $\calk$ can be characterized as
	\begin{align}
		\calk^* = \left\{
			\lambda =
			\left(\begin{matrix}
				\rho \\ s \\ z
			\end{matrix}\right)
			: \rho>0, 
				\left( z - \frac{1}{4\rho} T^*(s s\h) \right) \in \calc^*
					\right\}
				\cup
				\left\{ \lambda : \rho=0, s=0, z\in\calc^* \right\}.
	\end{align}
\end{lemma}
\begin{IEEEproof}
	See the appendix.
\end{IEEEproof}
It is clear that $\calk$ is not self-dual ($\calk\ne\calk^*$) and so
\eqref{conicproblem} is a non-symmetric conic program.

The cone $\calc$ and its dual are defined in terms of real-valued vectors
because this description simplifies the derivation of the method in Sec.
\ref{sec:method}. These sets are, however, more naturally understood from their
corresponding complex-valued forms. We therefore define the vector $u_\bbc =
(u_0, u_1+ju_{N}, u_2+ju_{N+1}, \ldots, u_{N-1}+ju_{2N-2})\T$ and use a
similar definition of $z_\bbc$.

The dual cone $\calc^*$ turns out to be the set of finite autocorrelation
sequences. An excellent introduction to this set and a number of
characterizations of it are given in \cite{alkire-autocorrelation} for the case
of real-valued sequences. Here we extend the definition to the complex-valued
case.
\begin{definition}
	A vector $z$ is a finite autocorrelation sequence if there exists a
	vector $q\in\bbc^N$ such that%
	\footnote{The complex conjugate of $q$ is denoted $\bar{q}$.}
	\begin{align}
		(z_\bbc)_k = \sum_{n=0}^{N-1-k} \bar{q}_n q_{n+k},
		\qquad k=0,\ldots,N-1.
		\label{acs_def}
	\end{align}
\end{definition}
In other words, $z$ is a finite autocorrelation sequence if
\begin{align}
	\ldots, 0, 0, (\bar{z}_\bbc)_{N-1}, (\bar{z}_\bbc)_{N-2}, \ldots, (\bar{z}_\bbc)_{1},
		(z_\bbc)_{0},
		(z_\bbc)_{1}, \ldots, (z_\bbc)_{N-1}, 0, 0, \ldots
		\label{acs}
\end{align}
is the autocorrelation sequence of some moving average process of order $N-1$
with filter coefficients $q_1, \ldots, q_{N-1}$ and input variance $|q_0|^2$.
It is well known from the theory of linear time-invariant systems, that if
\eqref{acs} is a valid autocorrelation, then it can be represented by a moving
average process (i.e., there exists a coefficient vector
$q$ such that \eqref{acs_def} holds).

A sequence is a valid autocorrelation
sequence if and only if its Fourier transform is non-negative \cite{KN:77}. The
Fourier transform of \eqref{acs} is
\begin{align}
	Z(\omega) = (z_\bbc)_0 + 2 \sum_{k=1}^{N-1} \re\!\left((z_\bbc)_k
	\exp(-j\omega k) \right),
\end{align}
for $\omega\in[0,2\pi)$.
Then $z\in\calc^*$ if and only if $Z(\omega)\ge0$ for all $\omega\in[0,2\pi)$.
The fast Fourier transform allows $Z(\omega)$ to be evaluated at a large number
of points on $[0,2\pi)$ in an efficient way.
Using Lemma \ref{Kstar_lemma} we therefore have a low-complexity method of
determining if $\lambda\in\calk^*$. This approach is approximate in the sense
that $Z(\omega)$ is sampled at a finite number of points on $[0,2\pi)$. The
approximation can be made arbitrarily accurate by increasing the number of
evaluated points. In our simulation study in Sec.~\ref{sec:numerical} it is
demonstrated that the approximation is of sufficient accuracy for our algorithm
to be utilized in practice.

We still haven't shown that the dual of the cone $\calc$ is indeed the set of
finite autocorrelation sequences. To that end, let $\tilde\calc$ be the set of
finite autocorrelation sequences and identify $u$ with $u_\bbc$.
Extending the approach of \cite{alkire-autocorrelation} to the complex-valued case, a vector $u$
is in the dual of $\tilde\calc$ if and only if $z\T u\ge0$ for every
$z\in\tilde\calc$, or, in other words, if and only if
\begin{align}
	z\T u &= \re(z_\bbc\h u_\bbc)
	= \re\!\left( \sum_{k=0}^{N-1}\sum_{n=0}^{N-1-k} (u_\bbc)_k q_n \bar{q}_{n+k}
	\right) = \frac{1}{2} q\T T(u) \bar{q} \ge 0
\end{align}
for every $q\in\bbc^N$. We can therefore identify $\calc$ with $\tilde\calc^*$.
Since $\tilde\calc$ is a proper cone, we have $\calc^* =
\tilde\calc^{**}=\tilde\calc$.

\subsection{Barrier Functions}
\label{sec:barriers}
IPMs are built on the idea of a barrier function
$F:\intr\calk\rightarrow\bbr$ associated to the cone $\calk$
($\intr\calk$ denotes the interior of $\calk$).
The barrier function must be a smooth and strictly convex%
\footnote{Hessian positive definite everywhere.}
function with $F(\mu_k)\rightarrow\infty$ for every sequence of points
$\mu_k\in\intr\calk$ with limit point $\tilde\mu\in\bdy\calk$, where
$\bdy\calk$ denotes the boundary of $\calk$.
The typical approach to IPMs also assumes that the barrier function
is logarithmically homogeneous (LH). $F$ is a LH barrier function
for the cone $\calk$ if there exists a $\theta_F>0$ such that $F(\alpha\mu) =
F(\mu) - \theta_F\log(\alpha)$ for all $\alpha>0$, $\mu\in\intr\calk$. The
value $\theta_F$ is called the degree of the barrier.

We will use the following well-known properties of a LH barrier function $F$
for $\calk$ \cite{BoVa:04, NeTo-primaldual, NeTo-selfscaled}:
If $\mu\in\intr\calk$, then
\begin{align}
	\left<-\nabla_\mu F(\mu), \mu \right> &= \theta_F,
	\label{graddotprod} \\
	-\nabla_\mu F(\mu) &\in\intr\calk^*,
	\label{gradindual}
\end{align}
where the gradient operator is defined as
$\nabla_\mu F = (\nabla_v F, \nabla_x F\T, \nabla_u F\T)\T$.
The gradient with respect to the complex vector $x=a+jb$ is to be understood
as%
\footnote{This is actually twice the Wirtinger derivative of $F$ with respect
to $\bar{x}$.}
$\nabla_x F = \nabla_af + j\nabla_bf$.

The usefulness of barrier functions is clear when considering their role in
path-following methods. A primal-only path-following method finds a solution to
\eqref{conicproblem} by iteratively solving
\begin{equation}
	\label{barrierproblem}
\begin{array}{ll}
	\text{minimize} & f(\mu) + t\ii F(\mu) \\[1mm]
	\text{subject to} & \mu\in\intr\calk
\end{array}
\end{equation}
for an increasing sequence of values $t>0$. In each step $\mu$ is initialized
with the solution of the previous step. This approach is desirable because each
step can be solved by an algorithm for unconstrained optimization such as
Newton's method.

In this paper we use the standard log-determinant barrier function for $\calk$:
\begin{align*}
	F(\mu) &=
		- \log
			\left| \left(\begin{matrix}
				T(u) & x \\ x\h & v
			\end{matrix}\right) \right|
			\shownumber\label{F}
			\\
		&= - \log|T(u)| - \log(v-x\h T\ii(u) x),
		\;\;\mathrm{for\ }\mu\in\intr\calk.
\end{align*}
It is easy to show that it is LH with degree $\theta_F=N+1$.

\subsection{Solvability}
We now consider conditions for the problem \eqref{conicproblem} to be solvable.
An optimization problem is solvable when a feasible point exists and its
objective is bounded below on the feasible set.
\begin{lemma}
	\label{bounded_lemma}
	The function $f(\mu)$ is bounded below on $\mu\in\calk$ if and only if
	$\tau=0$ or $\tau>0$ and $w\in\calc^*$.
\end{lemma}
\begin{IEEEproof}
	The case $\tau=0$ is trivial. Assume $\tau\ne0$ in the following.

	If $\tau<0$ or $w\notin\calc^*$ there exists $\mu\in\calk$ with $x=0$ such
	that $\tau v + \tau w\T u<0$. Then $\alpha\mu\in\calk$ for any
	$\alpha\ge0$ and $\lim_{\alpha\rightarrow\infty}f(\alpha\mu)=-\infty$,
	so $f(\mu)$ is unbounded below on $\mu\in\calk$.
	
	Conversely, if $\tau>0$ and $w\in\calc^*$, we have $\tau v\ge0$ and $\tau
	w\T u\ge0$
	for every $\mu\in\calk$. So $f(\mu)\ge0$ for $\mu\in\calk$.
\end{IEEEproof}
Since a primal feasible point always exists (take for example $v=1,x=0,u=e_0$),
the problem \eqref{conicproblem} is solvable if and only if $\tau=0$ or the
conditions in Lemma \ref{bounded_lemma} are fulfilled. These conditions can
easily be checked prior to executing the algorithm and we assume that the
problem is solvable in the following.

\subsection{Optimality Conditions}
With the conic modelling machinery in place we can begin to analyze the
solution of \eqref{problem} by considering the non-symmetric conic formulation
\eqref{conicproblem}.
The Lagrangian of \eqref{conicproblem} is
\begin{align}
	L(\mu,\lambda) = \norm{x-y}_2^2 + \tau(v+w\T u) - \left<\lambda,\mu\right>
\end{align}
and the dual is
\begin{equation}
	\label{dualproblem}
\begin{array}{ll}
	\text{maximize}
		& -\frac{1}{4}\norm{s}_2^2 - \re(y\h s) \\[1mm]
	\text{subject to}
		& \lambda\in\calk^*,\; \rho=\tau,\; z=\tau w.
\end{array}
\end{equation}
Notice that by taking the dual of \eqref{conicproblem} instead of
\eqref{problem}, the number of dual variables is reduced from $\calo(N^2)$ to
$\calo(N)$ (see \cite{Bhaskar:2013} for an explicit formulation of the dual of
\eqref{problem}). This is the reason why, from a computational point of view, it is
beneficial to work with the form \eqref{conicproblem} instead of
\eqref{problem}.

Since $f$ is convex, the Karush-Kuhn-Tucker (KKT) are necessary and sufficient
\cite[Sec. 5.9]{BoVa:04}
for variables $(\mu^\star, \lambda^\star)$ to be solutions of the primal and dual
problems \eqref{conicproblem} and \eqref{dualproblem}.
The KKT conditions are
\begin{align}
	\left\{ \begin{array}{l}
		\nabla_\mu L\!\left(\mu^\star,\lambda^\star\right) = 0 \\
		\mu^\star \in\calk \\
		\lambda^\star \in\calk^* \\
		\left<\lambda^\star,\mu^\star\right> = 0
	\end{array} \right\}.
	\label{kkt}
\end{align}

\newcommand{\ttt}{^{(t)}}
Instead of directly solving the KKT conditions, our primal-dual IPM finds solutions $(\mu\ttt, \lambda\ttt)$ of the augmented KKT conditions
\cite{BoVa:04, NeTo-ipm}
\begin{align}
	\left\{ \begin{array}{l}
		\nabla_{\mu}L\!\left(\mu\ttt,\lambda\ttt\right) = 0 \\
		\mu\ttt \in \intr\calk \\
		\lambda\ttt \in \intr\calk^* \\
		\lambda\ttt = - t\ii \nabla_{\mu} F\!\left(\mu\ttt\right)
	\end{array} \right\}
	\label{akkt}
\end{align}
for an increasing sequence of values $t>0$. It is easy to realize that
$\left(\mu\ttt,\lambda\ttt\right)$ solves \eqref{akkt} only if $\mu\ttt$
is a solution of the barrier problem \eqref{barrierproblem}. This observation
provides the link between primal-only barrier methods and primal-dual IPMs. The
set of values $\left\{\left( \mu\ttt,\lambda\ttt \right) : t>0\right\}$ is
known as the primal-dual central path. The primal-dual central path converges
to the desired solution in the sense that $\lim_{t\rightarrow\infty}\left(
\mu\ttt,\lambda\ttt \right) = \left( \mu^\star, \lambda^\star \right)$
\cite{BoVa:04, NeTo-ipm}.

The last condition in \eqref{akkt} is known as the augmented complementary
slackness condition. It follows from \eqref{gradindual} that the
second and fourth condition in \eqref{akkt} together imply
$\lambda\ttt\in\intr\calk^*$, so the third condition can be dropped.

From \eqref{graddotprod} it follows that the duality gap for the primal-dual
problems \eqref{conicproblem} and \eqref{dualproblem} at a point on the
primal-dual central
path is $\left<\lambda\ttt,\mu\ttt\right>=\theta_F/t$
\cite{BoVa:04,NeTo-selfscaled}. So solving the augmented KKT gives a
primal feasible solution $\mu\ttt$ which is no more than $(N+1)/t$ suboptimal.
Consequently, an arbitrarily accurate solution can be obtained by solving
\eqref{akkt} for sufficiently large $t$.

\subsection{Obtaining a Solution of the Augmented KKT Conditions}
We now define $v\ttt, x\ttt, u\ttt, \rho\ttt, s\ttt$ and  $z\ttt$ as the
entries of $\mu\ttt$ and $\lambda\ttt$.
By solving the first equation in \eqref{akkt} (the stationarity condition) we get
\begin{align}
	\rho\ttt=\tau,\quad z\ttt=\tau w,\quad s\ttt=2(x\ttt-y).
	\label{lambda_t}
\end{align}
We continue by writing out the last condition of \eqref{akkt}. Solve for $v\ttt$ and $x\ttt$
and insert the relations above to get
\begin{align}
	v\ttt &= (\tau t)\ii + \left(x\ttt\right)\h T\ii\!\left(u\ttt\right) x\ttt
	\label{v_t} \\
	x\ttt &= T\!\left(u\ttt\right) T\ii\!\left(u\ttt+2\ii\tau e_0\right) y.
	\label{x_t}
\end{align}
Finally, solve $z\ttt=-t\ii\nabla_uF\!\left(\mu\ttt\right)$
for $u\ttt$ and insert the above to obtain
\begin{align}
	\label{compslack_u}
	\tau w
	- \tau T^*\!\left( \phi\phi\h \right)
	- t\ii T^*\!\left( T\ii(u\ttt) \right) = 0,
\end{align}
where $\phi = T\ii\!\left(u\ttt+2\ii\tau e_0\right) y$.

For a given $t>0$ the corresponding point on the primal-dual central path can
be obtained as follows: First a solution $u\ttt$ of \eqref{compslack_u} that
fulfills $u\ttt\in\intr\calc$ is found (existence of such a solution is shown
below).  Then the point $\left(\mu\ttt,\lambda\ttt\right)$ is obtained by
inserting into \eqref{lambda_t}, \eqref{v_t} and \eqref{x_t}. It is easy to
show from $u\ttt\in\intr\calc$ that $\mu\ttt\in\intr\calk$ and so
$\left(\mu\ttt,\lambda\ttt\right)$ solves \eqref{akkt} and it is a primal-dual
central point.

How can we obtain a solution $u\ttt\in\intr\calc$ of \eqref{compslack_u}?  The
left-hand side of \eqref{compslack_u} is recognized as the gradient of $h_t(u)
= g(u) + t\ii G(u)$, with
\begin{align}
	g(u) &= \tau w\T u + \tau y\h T\ii\!\left(u+2\ii\tau e_0\right) y
	\label{g} \\
	G(u) &= - \log|T(u)|.
	\label{G}
\end{align}
Now consider the barrier problem
\begin{equation}
	\label{min_g}
\begin{array}{ll}
	\text{minimize} & h_t(u) \\[1mm]
	\text{subject to} & u\in\intr\calc.
\end{array}
\end{equation}
The gradient of $h_t$ vanishes at the solution of \eqref{min_g} because $G$ is
a LH barrier function for $\calc$.
So solving \eqref{compslack_u} with $u\ttt\in\intr\calc$ is equivalent to
solving \eqref{min_g}.
Since we have assumed that the problem \eqref{conicproblem} is solvable, then so is
\eqref{min_g} (thus proving that there exists a $u\ttt\in\intr\calc$ that
solves \eqref{compslack_u}).

The idea of the primal-dual IPM presented in the following section is to use an
iterative algorithm for unconstrained optimization (either Newton's method or a
quasi-Newton method) to solve \eqref{min_g}. However, we do not need to exactly
solve \eqref{min_g} for a sequence of values $t>0$. In each iteration of the
solver the value of $t$ can be updated in a dynamic manner based on the duality
gap.

\section{The Primal-Dual Interior-Point Method}
\label{sec:method}

\newcommand{\pre}{_{i-1}}
\newcommand{\cur}{_{i}}
\newcommand{\nex}{_{i+1}}

We now outline FastAST, a primal-dual IPM for the solution of
\eqref{conicproblem}. Let $(\mu\cur,u\cur,\lambda\cur,t\cur)$ denote
$(\mu,u,\lambda,t)$ in iteration $i$.  The proposed method is given in
Algorithm \ref{alg}. In the remainder of this section, each step of the
algorithm is discussed in detail.

\begin{figure}[t]
\removelatexerror
\begin{algorithm}[H]
	\DontPrintSemicolon
	\KwParam{$\gamma>1$.}
	\KwIn{
		Initial values $u_0\in\calc$ and $t_1>0$.
	}
	Set objective lower bound $f_|LB|=-\infty$. \;
	\For{$i=1,2,\ldots$}{
		Determine the search direction $\Delta u$. \;
		Perform a line search along $\Delta u$ to obtain the step size
			$\alpha$. \;
		Update estimate $u\cur=u\pre + \alpha\Delta u$. \;
		Form primal-dual variables $(\mu\cur, \lambda\cur)$ using \eqref{lambda_t}, \eqref{v_t} and
			\eqref{x_t}. \;
		\If{$\lambda\cur\in\calk^*$}{
			Update lower bound on objective
				$f_|LB| = \max\!\left(f_|LB|, -\frac{1}{4}\norm{s\cur}_2^2
					- \re(y\h s\cur)\right)$. \;
		}
		Evaluate duality gap $\eta\cur = f(\mu\cur) - f_|LB|$. \;
		Terminate if the stopping criterion is satisfied. \;
		Update barrier parameter $t\nex = \max\!\left( t\cur, \gamma \frac{N+1}{\eta\cur}
				\right)$. \;
	}
	\KwOut{
		Primal-dual solution $(\mu\cur, \lambda\cur)$.
	}
	\caption{Primal-dual IPM for fast atomic norm soft
	thresholding (FastAST).}
	\label{alg}
\end{algorithm}
\end{figure}

Low-complexity evaluation of the steps in FastAST are presented in
Sec.~\ref{sec:complexity}. With these approaches, the computational complexity
is dominated by the evaluation of the search direction. For this step we
propose to use either Newton's method or a quasi-Newton method. The
quasi-Newton method has much lower computational complexity per iteration and
is also faster in practice. It is, however, not able to obtain a solution of
high accuracy. If a highly accurate solution is required, Newton's method is
preferred. We refer to the numerical evaluation in Sec.~\ref{sec:numerical} for
a detailed discussion thereof.

\subsection{Determining the Search Direction Using Newton's Method}
Applying Newton's method to solve \eqref{min_g} we get the search direction
\begin{align}
	\Delta u = - \left( \nabla^2_u h_{t\cur}(u\pre) \right)\ii
		\nabla_u h_{t\cur}(u\pre),
	\label{du_newton}
\end{align}
where $\nabla^2_u h_{t\cur}(u\pre)$ denotes the Hessian of $h_{t\cur}$
evaluated at $u\pre$. As discussed in Sec.~\ref{sec:complexity} the Hessian
can be evaluated in $\calo(N^3)$ flops and the same cost is required for
solution of the system \eqref{du_newton}.

\subsection{Determining the Search Direction Using L-BFGS}
\newcommand{\kkk}{_{k}}

In scenarios with large $N$ the computation time for evaluation of the Newton
search direction can become prohibitively large. In these cases we propose to
use the limited-memory Broyden–Fletcher–Goldfarb–Shanno (L-BFGS) algorithm
\cite{nocedal-lbfgs} for the solution of \eqref{min_g}. L-BFGS posses two key
properties that are instrumental in obtaining an algorithm with low
per-iteration computational complexity: \textit{1)} it uses only gradient
information and the gradient of $h_{t\cur}$ can be evaluated with low
computational complexity;%
\footnote{To speed up convergence, our implementation also uses an
approximation of the diagonal of the Hessian of $h_{t\cur}$.}
and \textit{2)} by
appropriately modifying the L-BFGS two-loop recursion, it can be used for the
solution of \eqref{min_g} in a computational efficient manner even though $t$
is increased in every iteration. It is this property, and not the limited
memory requirements, that makes L-BFGS preferable over other quasi-Newton
methods (such as vanilla BFGS) for our purposes.

In relation to the second property, note that since $t\cur\neq t\pre \neq
\dots$, the normal formulation of L-BFGS does not apply. A simple modification
of the L-BFGS two-loop recursion \cite{nocedal-lbfgs} overcomes this
limitation. At the end of the $i$th iteration, the following difference vectors
are calculated and saved for later use:
\begin{align}
	r\cur &= u\cur - u\pre \label{r_i} \\
	q\cur &= \nabla_u g(u\cur) - \nabla_u g(u\pre) \\
	Q\cur &= \nabla_u G(u\cur) - \nabla_u G(u\pre) \label{qG_i}.
\end{align}
This set of vectors is retained for $M$ iterations.  The modified two-loop
recursion in Algorithm \ref{lbfgs} can then be used to calculate the search
direction $\Delta u$. This algorithm calculates the normal L-BFGS search
direction for minimization of $h_{t\cur}$, as if $t\cur = t\pre = \ldots$.
That can be achieved because L-BFGS only depends on $t\cur$ through
$\nabla_uh_{t\cur}(u_k)=\nabla_ug(u_k) + t\cur\ii \nabla_uG(u_k)$, for
$k=i-1,\ldots,i-M-1$. The gradients $\nabla_ug(u_k)$ and
$\nabla_uG(u_k)$ need only be calculated once to allow $\nabla_uh_{t\cur}(u_k)$ to be
calculated for any value of $t\cur$.

\begin{figure}[t]
\removelatexerror
\begin{algorithm}[H]
	\DontPrintSemicolon
	\KwParam{Number of saved difference vectors $M$.}
	\KwIn{Current iteration number $i$ and parameter $t\cur$.
		Saved difference vectors $r_k, q_k, Q_k$ for $k=i-1,i-2,\ldots,\max(i-M,1)$.
		Current gradient vector $\nabla_u h_{t\cur}(u\pre)$ and initial
		Hessian approximation $\hat H\cur$.
	}
	$d \gets - \nabla_u h_{t\cur}(u\pre)$ \;
	\For{$k=i-1, i-2, \ldots, \max(i-M,1)$}{
		$\psi\kkk \gets q\kkk + t\cur\ii Q\kkk $ \;
		$\sigma\kkk \gets \frac{r\kkk\T d}{r\kkk\T \psi\kkk}$ \;
		$d \gets d - \sigma\kkk \psi\kkk$
	}
	$d \gets \hat H\cur\ii d$ \;
	\For{$k=\max(i-M,1), \max(i-M,1)+1, \ldots, i-1$}{
		$\beta\kkk \gets \frac{\psi\kkk\T d}{\psi\kkk\T r\kkk}$ \;
		$d \gets d + r\kkk \left( \sigma\kkk - \beta\kkk \right)$
	}
	\KwOut{
		Search direction $\Delta u = d$. \;
	}
	\caption{Modified L-BFGS two-loop recursion for calculation of the search
	direction.}
	\label{lbfgs}
\end{algorithm}
\end{figure}


In each iteration the initial Hessian $\hat H\cur$ should be chosen as an
approximation of the Hessian of $h_{t\cur}$ evaluated at $u\pre$. It is the matrix upon
which L-BFGS successively applies rank-2 updates to form the Hessian
approximation that is used for calculating the search direction.
An easy, and popular, choice for the initial Hessian is the identity matrix
$\hat H\cur=I$.
Through numerical experiments we have seen that this choice leads to slow
convergence. It turns out that the slow convergence is caused by
the scaling of the Hessian, leading to non-acceptance of a full Newton step
(i.e., $\alpha$ is selected much smaller than 1). Using a diagonal approximation of the true
Hessian remedies this, but, unfortunately, it cannot be calculated with low
computational complexity. (Our best attempt at devising a fast
evaluation of the Hessian diagonal yielded cubic complexity $\calo(N^3)$, the
same as evaluation of the full Hessian.) Instead our algorithm uses the
following heuristic approximation of the diagonal Hessian
\begin{align}
	\label{Hi}
	\hat H\cur = \diag\!\left( 1, \frac{N-1}{2N}, \ldots, \frac{1}{2N},
		\frac{N-1}{2N}, \ldots, \frac{1}{2N} \right)
			\left(\nabla^2_u h_{t\cur}(u\pre)\right)_{0,0},
\end{align}
where $\left(\nabla^2_u h_{t\cur}(u\pre)\right)_{0,0}$ is the $(0,0)$th entry
of the true Hessian evaluated at $u\pre$.
This approximation can be calculated with low computational complexity as
demonstrated in Sec.~\ref{sec:complexity}. The approximation is motivated as
follows: The diagonal entries are scaled according to the
number of times the corresponding entry of $u$ appears in $T(u)$. This scaling
resembles that in the biased autocorrelation estimate (except for a factor of
$2$ caused by the scaling of the diagonal in the definition of $T(u)$). In our
numerical experiments, we have observed the above approximation to be fairly
accurate; each entry typically takes a value within $\pm50\,\%$ of the true
value. To this end we note that only a crude approximation is needed, since the
role of $\hat H\cur$ is to account for the scaling of the problem. Our numerical
investigation suggests that using the approximation \eqref{Hi} leads to only
marginally slower convergence, compared to using a diagonal Hessian
approximation using the diagonal of the true Hessian.

A final note on our adaptation of L-BFGS is that the usual observations
regarding positive definiteness of the approximated Hessian remain valid.
First note that the objective upon which L-BFGS is applied ($h_{t\cur}$) is
a strictly convex function for $u\in\intr\calc$. It follows that 
the initial Hessian approximation $\hat H\cur$ is positive definite. Also, the
curvature condition $r_k\T\psi_k>0$ is valid for all $k$. Then the
approximated Hessian is positive definite and the calculated search direction
$\Delta u$ is a descent direction \cite{nocedal-lbfgs, Nocedal99}.

\subsection{Line Search}
The line search along the search direction $\Delta u$ is a simple backtracking
line search starting at $\alpha=1$. A step size is accepted if the new point is
strictly feasible, i.e., if $u\pre + \alpha\Delta u\in\intr\calc$. It is then easy to
show that the primal solution $\mu\cur$ calculated from inserting $u\cur$ into
\eqref{v_t} and \eqref{x_t} is strictly primal feasible ($\mu\cur\in\intr\calk$).

To guarantee that the objective is sufficiently decreased,
the Armijo rule is also required for acceptance of a step size $\alpha$:
\begin{align}
	h_{t\cur}(u\pre+\alpha\Delta u)
		&\le h_{t\cur}(u\pre) + c \alpha \Delta u\T \nabla_u h_{t\cur}(u\pre),
\end{align}
where $0<c<1$ is a suitably chosen constant.

\subsection{The Duality Gap and Update of \texorpdfstring{$t$}{t}}
The line search guarantees that the primal solution is strictly feasible in all
iterations, i.e., that $\mu\cur\in\intr\calk$. Dual feasibility of a solution $\lambda\cur$
obtained from \eqref{lambda_t} is not guaranteed. The algorithm
therefore checks for $\lambda\cur\in\calk^*$ using the approximate approach
described in Sec. \ref{dualmembership}. 

Let $f^\star$ denote the optimal value of the problem \eqref{conicproblem}.
If $\lambda\cur$ is dual feasible, the objective of the dual
\eqref{dualproblem} provides a lower bound on the optimal value, i.e.,
\begin{align}
	f^\star \ge -\frac{1}{4}\norm{s\cur}_2^2 - \re(y\h s\cur).
\end{align}
The algorithm always retains the largest lower bound it has encountered in $f_|LB|$.
From the lower bound, a duality gap $\eta\cur$ can be evaluated in each iteration:
\begin{align}
	\eta\cur = f(\mu\cur) - f_|LB|.
	\label{eta}
\end{align}
This value gives an upper bound on the sub optimality of the
solution $\mu\cur$, i.e., $f(\mu\cur)-f^\star\le\eta\cur$.

Recall that the algorithm is ``aiming'' for a solution of the augmented KKT
conditions \eqref{akkt}. At this
solution, the duality gap is $\theta_F/t\nex$. The next value of $t$ can then
be determined so that the algorithm is aiming for a suitable (not too
large, not too small) decrease in the duality gap, i.e., we select $t\nex$ such
that $\eta_i/\gamma = \theta_F/t\nex$ for some preselected $\gamma>1$. To
guarantee convergence it is also imposed that $t\cur$ is a non-decreasing
sequence.

\subsection{Termination}
The duality gap provides a natural stopping criterion. The proposed algorithm
terminates based on either the duality gap ($\eta\cur<\varepsilon_|abs|$) or
the relative duality gap ($\eta\cur / f(\mu\cur) < \varepsilon_|rel|$).  The
relative duality gap is a sensible stopping criterion because $f(\mu)\ge0$ as
is seen in the proof of Lemma \ref{bounded_lemma}.

Algorithm \ref{alg} is guaranteed to terminate at a point that fulfills either
of the two stopping criteria listed above.
To see why that is the case, consider a scenario where $t\cur$ converges to some
finite constant $\tilde t$ as $i\rightarrow\infty$. Then, as $i\rightarrow\infty$,
the algorithm implements L-BFGS with a backtracking line search to minimize
$h_{\tilde t}$. Thus $u\cur$ converges to the minimizer $u^{(\tilde t)}$ of
$h_{\tilde t}$. Let $(\mu^{(\tilde t)}, \lambda^{(\tilde t)})$ denote the
corresponding primal and dual variables calculated from \eqref{v_t},
\eqref{x_t} and \eqref{lambda_t}.

Now, $(\mu^{(\tilde t)}, \lambda^{(\tilde t)})$ constitute a solution to \eqref{akkt}
with $t=\tilde t$. Then $\lambda^{(\tilde t)}\in\intr\calk^*$ follows from
\eqref{gradindual}. 
Further, we have from \eqref{eta}, \eqref{lambda_t} and \eqref{graddotprod} that the duality
gap $\eta\cur$ converges to $\left< \mu^{(\tilde t)}, \lambda^{(\tilde t)}
\right> = \theta_F / \tilde t$ as $i\rightarrow\infty$. However, that
implies $t\nex=\gamma\theta_F/\eta_i=\gamma\tilde t>\tilde t$ in the limit, a
contradiction to the assumption that $t\cur$ converges to $\tilde t$. This
means that $t\cur$ does not converge to a finite value and, as it is
non-decreasing, it must diverge to $+\infty$. It is also evident that the
duality gap $\eta\cur\rightarrow0$ as $t\cur\rightarrow\infty$, and so either
of the stopping criteria are eventually fulfilled.

\subsection{Initialization}
FastAST must be initialized with a primal variable $u_0\in\calc$ and
a barrier parameter $t_1>0$.
To determine a suitable value of the initial barrier parameter $t_1$ we first identify a
primal-dual feasible point from which the duality gap can be evaluated. A
primal-dual feasible point can be obtained by assuming%
\footnote{The problem \eqref{problem} is solvable if and only if $w\in\calc^*$.
The restriction to the interior has no practical effect.}
$w\in\intr\calc^*$ and iterating these steps:
\begin{enumerate}
	\item Set $u=10\lVert y\rVert_2^2/N, 0, \ldots, 0)\T$.
	\item Calculate $(\mu,\lambda)$ from $u$ based on \eqref{lambda_t},
		\eqref{v_t} and \eqref{x_t}.
	\item If $\lambda\in\calk^*$, terminate, otherwise double the first entry
		of $u$ and go to step 2.
\end{enumerate}
The value of $u$ in Step 1 has been chosen heuristically.
It is easy to see that $u$ stays primal feasible throughout.
It is guaranteed that a dual feasible point is reached because
$u\rightarrow(\infty,0,\ldots,0)\T$. Then, following \eqref{x_t},
we have $x\rightarrow y$ and so $s\rightarrow 0$. Considering the result in
Lemma \ref{Kstar_lemma} and the assumption $w\in\intr\calc^*$ we get that
$\lambda$ converges to a point $\tilde\lambda\in\intr\calk^*$.

When a primal-dual feasible point $(\mu,\lambda)$ has been found the
corresponding duality gap is $\eta_0=\left< \mu, \lambda \right>$.  The initial
value of the barrier parameter is selected as $t_1 = \gamma\theta_F / \eta_0$.
The corresponding value of $u$ is used as the initial value of the primal
variable $u_0$.

\section{Fast Computations}
\label{sec:complexity}
For brevity iteration indices are dropped in the following.
The computationally demanding steps of FastAST (Alg.~\ref{alg}) all involve the
determinant or the inverse of Toeplitz matrices $T(u)$ and $T(u+2\ii\tau e_0)$.
In this section we demonstrate how the Toeplitz structure can be exploited to
significantly reduce the computational complexity of these evaluations. The
use of such structure for fast solution of optimization problems
have previously been seen \cite{HansenFR17, musicus-fastmlm}, including for
evaluation of the gradient and Hessian of the barrier function $G$ 
\cite{genin-optimization, alkire-autocorrelation}.


\subsection{Fast Algorithms for Factorizing a Toeplitz Inverse}
Our computational approach is based on the following factorizations of Toeplitz
inverses.  The Gohberg-Semencul formula \cite{ammar-generalized,
gohberg-convolution} gives a factorization of the inverse of a Toeplitz matrix
$T(u)$,
\begin{align}
	T\ii(u) = \delta_{N-1}\ii (U\h U - VV\h),
	\label{gohberg-semencul}
\end{align}
where the entries of Toeplitz matrices $U$ and $V$ are
\begin{align}
	U_{n,m} &= \rho_{N-1+n-m}, \\
	V_{n,m} &= \rho_{n-m-1},
\end{align}
for $n,m=0,\ldots,N-1$. Note that $\rho_n=0$ for $n<0$ and $n>N-1$; thus $U$ is
unit upper triangular ($\rho_{N-1}=1$) and $V$ is strictly lower triangular.

The values $\delta_n$ and $\rho_n$ for $n=0,\ldots,N-1$ can be computed with a
generalized Schur algorithm in $\calo(N\log^2N)$ flops \cite{ammar-generalized}.
Alternatively, the Levinson-Durbin algorithm can be used to obtain the
decomposition in $\calo(N^2)$ flops. The latter algorithm is significantly
simpler to implement and is faster for small $N$. In \cite{ammar-numerical} it
is concluded that the Levinson-Durbin algorithm requires fewer total operations
than the generalized Schur algorithm for $N\le256$.

We will also use a Cholesky factorization of $T\ii(u)$,
\begin{align}
	T\ii(u)=PDP\h\,
	\label{PDP}
\end{align}
where $P$ is unit upper triangular and $D$ is diagonal.
The matrix $D=\diag(\delta_0\ii,\ldots,\delta_{N-1}\ii)$ is inherently computed
when the generalized Schur algorithm is executed \cite{ammar-generalized}. The
generalized Schur algorithm does not compute the matrix $P$. The
Levinson-Durbin algorithm inherently computes both $P$ and $D$, a property
which we exploit for evaluation of the Hessian of the barrier function $G$.

In the following we let $\rho_0, \ldots, \rho_{N-1}$ and $\delta_0, \ldots,
\delta_{N-1}$ be the entries obtained by executing the generalized Schur or
Levinson-Durbin algorithm with either $T\ii(u)$ or $T\ii(u+2\ii\tau e_0)$;
which one will be clear from the context.

\subsection{Evaluating the Objective and the Primal Variables}

We first discuss evaluation of the objective $h_{t}(u) = g(u) + t\ii G(u)$.
Since $P$ in \eqref{PDP} has unit diagonals, it is easy to obtain
\begin{align}
	G(u) = - \log|T(u)| = - \sum_{n=0}^{N-1} \log\delta_n.
\end{align}

To evaluate $g(u)$ insert \eqref{gohberg-semencul} into \eqref{g} and realize
that all matrix-vector products involve Toeplitz matrices. Vector multiplication onto
a Toeplitz matrix can be performed using the fast Fourier transform (FFT) in
$\calo(N\log N)$ flops (such products are convolutions,
see e.g.~\cite{alkire-autocorrelation} for details). In conclusion, the
dominant cost of evaluating $h_t(u)$ is the execution of the generalized Schur
(or Levinson-Durbin) algorithm.

Evaluating the primal variables $v^{(t)}$ and $x^{(t)}$ in
\eqref{v_t} -- \eqref{x_t} similarly amounts to
vector products onto Toeplitz matrices.

The line search in Algorithm \ref{alg} must check for $u\in\calc$, i.e., if
$T(u)\succ0$. The generalized Schur (or Levinson-Durbin) algorithm can again be
used here, as $T(u)\succ0$ if and only if $\delta_n>0$ for $n=0,\ldots,N-1$.

\subsection{Evaluating the Gradients}
The following gradients must be evaluated in each iteration of Algorithm
\ref{alg}:
\begin{align}
	\nabla_u g(u) &=
		\tau w
		- \tau T^*\!\left( \phi\phi\h \right) \\
	\nabla_u G(u) &=
		- T^*\!\left( T\ii(u) \right).
\end{align}

We first consider the term $T^*(\phi\phi\h)$. The vector $\phi$ can be
evaluated with low complexity (confer the evaluation of primal variables,
above). Let $\beta_n\in\bbc$ denote the sum over the $n$th upper diagonal of
$\phi\phi\h$ for $n=0,\ldots,N-1$, i.e.,
\begin{align}
	\beta_n &=
		\sum_{m=0}^{N-1-n} (\phi\phi\h)_{m,m+n}
		= \sum_{m=0}^{N-1-n} \phi_m\bar\phi_{m+n}.
		\label{sum_T_diag}
\end{align}
It is recognized that the values $\beta_0,\ldots,\beta_{N-1}$ can be calculated as a
correlation, which can be implemented using FFTs in
$\calo(N\log N)$ flops.
Then $T^*(\phi\phi\h)$ can be obtained by concatenating and
scaling the real and imaginary parts of $\beta$,
\begin{align}
	T^*(\phi\phi\h) =
	(2\beta_0, 2\re(\beta_1), \dots, 2\re(\beta_{N-1}),
	2\im(\beta_1), \dots, 2\im(\beta_{N-1}))\T\,.
	\label{a-to-T}
\end{align}

Now consider evaluation of the term $T^*\!\left( T\ii(u) \right)$.
The sum over the $n$th upper diagonal of $T\ii(u)$ is denoted as $\tilde\beta_n$ and
can be rewritten as
\begin{align}
	\tilde\beta_n &= \sum_{m=0}^{N-1-n} (T\ii(u))_{m,m+n}
		=\delta_{N-1}\ii \sum_{k=0}^{N-1} 
	(n-N + 2(k+1)) \rho_{k} \bar\rho_{k + n}\,,
\end{align}
see \cite{HansenFR17,musicus-fastmlm} for details. The above is recognized as
two correlations, thus allowing a low-complexity evaluation. The vector
$T^*\!\left( T\ii(u) \right)$ is found by concatenating and scaling the real
and imaginary parts of $\tilde\beta$, analogously to \eqref{a-to-T}.

\subsection{Evaluating the Full Hessian}
When Newton's method is used to determine the search direction, the Hessian of
$h_t$ must be evaluated. We now derive an approach to calculate the Hessians
of $g$ and $G$, from which the required Hessian is easily found.

The $(n,m)$th entry of the Hessian of $g$ is
\begin{align}
	\label{hess_g}
	\left( \nabla^2_ug(u) \right)_{n,m}
		= 2\tau \phi\h (E_n+E_n\h)
			T\ii(u+2\ii\tau e_0)
			(E_m+E_m\h) \phi,
\end{align}
where
\begin{align}
	E_n =
	\begin{cases}
		I			& n = 0 \\
		\tilde E^n			& 1 \le n \le N-1 \\
		-jE_{n-N+1}	& N \le n \le 2N-1.
	\end{cases}
\end{align}
The matrix $\tilde E$ is the lower shift matrix, i.e., it has ones on the lower
subdiagonal and zeros elsewhere. Note that $T(e_n)=E_n+E_n\h$. The $m$th column
of the Hessian is then
\begin{align}
	\left( \nabla^2_ug(u) \right)_{m}
		&= \tau
		T^*\! \left( d_m \phi\h + \phi d_m\h  \right),
\end{align}
where we let $d_m$ denote a vector $d_m = T\ii(u+2\ii\tau e_0) (E_m+E_m\h)
\phi$. Assuming the decomposition \eqref{gohberg-semencul} is available, a
column of the Hessian can be calculated in $\calo(N\log N)$ flops by explicitly
forming $d_m$ and performing sums over diagonals (as in \eqref{sum_T_diag}).
The full Hessian of $g$ is then obtained in $\calo(N^2\log N)$ flops.

To evaluate the Hessian of the barrier function $G$ we generalize the
approach of \cite{alkire-autocorrelation} to the complex-valued case. The
$(n,m)$th entry of the Hessian is
\begin{align}
	\label{hess_G}
	\left( \nabla_u^2 G(u) \right)_{n,m}
		&= \tr\!\left( T\ii(u) (E_n+E_n\h)
			T\ii(u) (E_m+E_m\h) \right).
\end{align}
Define the $N\times N$ matrices $A$ and $B$ with entries
\begin{align}
	A_{n,m} &= 2 \tr\!\left( T\ii(u) E_n T\ii(u) E_m \right)\\
	B_{n,m} &= 2 \tr\!\left( T\ii(u) E_n T\ii(u) E_m\T \right).
\end{align}
Then the Hessian can be written in the form
\begin{align}
	\label{hess_G_block}
	\nabla^2_uG(u) =
	\left(\begin{matrix}
		\re(A+B)		&	\re(-jAJ\T) \\
		\re(-jJA-jJB)	&	\re(-JAJ\T + JBJ\T)
	\end{matrix}\right),
\end{align}
where $J$ is a matrix that removes the first row, i.e., $J=(0,I)$, where $0$
is a column of zeros and $I$ is the $(N-1)\times(N-1)$ identity matrix.

At this point, we need a fast way of evaluating matrices $A$ and $B$.
Define the discrete Fourier transform matrix $W\in\bbc^{N_|FFT|\times N}$ with
entries
\begin{align}
	W_{n,m} = \exp(-j2\pi nm / N_|FFT|),
\end{align}
where $N_|FFT|$ is chosen such that $N_|FFT|\ge2N-1$.
Recall that the Levinson-Durbin algorithm gives the decomposition $T\ii(u) =
PDP\h$, from which $T\ii(u) = RR\h$ is obtained by calculating $R =
PD^{\frac{1}{2}}$.
Let $S_n$ denote the discrete Fourier transform of the $n$th column of $R$ (denote this
column $R_n$), i.e., $S_n = WR_n$.
Then by straight-forward generalization of the derivation in
\cite{alkire-autocorrelation} to the complex-valued case, we get that $A$ and
$B$ can be written in the forms
\begin{align}
	A &= \frac{2}{N_|FFT|^2} W\T \left(
		\left( \sum_{l=0}^{N-1} S_lS_l\h \right) \odot
		\left( \sum_{l=0}^{N-1} S_lS_l\h \right)
	\right) W
	\label{A} \\
	B &= \frac{2}{N_|FFT|^2} W\T \left(
		\left( \sum_{l=0}^{N-1} S_lS_l\h \right) \odot
		\left( \sum_{l=0}^{N-1} S_lS_l\h \right)
	\right) \bar W,
	\label{B}
\end{align}
with $\odot$ denoting the Hadamard (entrywise) product. Using \eqref{A} --
\eqref{B} the Hessian of $G$ can be evaluated in $\calo(N^3)$ flops.

\subsection{Evaluating the Diagonal Hessian Approximation}
The L-BFGS variant of FastAST uses the approximation of the Hessian diagonal \eqref{Hi}
which requires calculation of the first entry of the Hessian
\begin{align}
	\left(\nabla^2_u h_{t}(u)\right)_{0,0}
		=
	\left(\nabla^2_u g(u)\right)_{0,0}
	+ \frac{1}{t} \left(\nabla^2_u G(u)\right)_{0,0}.
\end{align}
An $\calo(N\log N)$ evaluation of the first term is easily obtained from
\eqref{hess_g}. The second term can be evaluated based on \eqref{hess_G_block}, but
a more efficient way is as follows: From \eqref{hess_G}
we have
\begin{align}
	\left( \nabla_u^2 G(u) \right)_{0,0}
		&=
		4 \tr\!\left( T\ii(u) T\ii(u) \right).
		\label{hess_G_first}
\end{align}
The matrix $T\ii(u)$ can be formed explicitly in $\calo(N^2)$ flops using 
the Trench algorithm \cite{trench-algorithm, golub-matrix}. However, since the
decomposition \eqref{gohberg-semencul} is already available in our setting it
is much easier to form $T\ii(u)$ from it by writing for $n=0,\ldots,N-1$ and
$m=0,\ldots,N-1-n$
{\small
\begin{align*}
	T\ii(u)_{m,m+n}
		= \delta_{N-1}\ii
		\left( \sum_{k=0}^m
			\bar\rho_{N-1-k} \rho_{N-1-(k+n)}
			- \rho_{k-1} \bar\rho_{k+n-1}
		\right),
\end{align*}
}
i.e., $T\ii(u)$ is ``formed along the diagonals''. By implementing the above
sum as a cumulative sum, the complete matrix $T\ii(u)$ is formed in
$\calo(N^2)$ flops.
Note that since $T(u)$ is both Hermitian and persymmetric, then so is
$T\ii(u)$. This means that only one ``wedge'' of the matrix, about
$N/4$ entries, must be calculated explicitly \cite{golub-matrix}.

The trace in \eqref{hess_G_first} is evaluated by taking the magnitude
square of all entries in $T\ii(u)$ and summing them.

\subsection{Analysis of Computational Complexity}

To summarize the computational complexity of an implementation of Alg.
\ref{alg} based on the low-complexity evaluations above, consider each of the two
variants for determining the search direction.
\begin{itemize}
	\item FastAST Newton: The computation time is asymptotically dominated by
		evaluation and inversion of the Hessian, i.e., $\calo(N^3)$
		flops.
	\item FastAST L-BFGS: The computation time is asymptotically dominated by the
		$\calo(MN)$ modified L-BFGS two-loop recursion in Alg. \ref{lbfgs} or
		by the $\calo(N^2)$ evaluation of the diagonal Hessian
		approximation.
\end{itemize}
When using the Newton search direction, the decomposition \eqref{PDP} is required and the
Levinson-Durbin algorithm must therefore be used to evaluate the factorization
of the Toeplitz inverse. When using the L-BFGS search direction either the
generalized Schur or the Levinson-Durbin algorithm can be used. The choice does
not affect the asymptotic computational complexity, but one may be faster than
the other in practice.

\section{Numerical Experiments}
\label{sec:numerical}
\subsection{Setup \& Algorithms}

\begin{table}[tbp]
	\centering
	\small
	\begin{tabular}{lcc}
		\toprule
		Variant & L-BFGS & Newton \\
		\midrule
		Number of saved difference vectors $M$
			& $2N-1$	& - \\
		Armijo parameter $c$
			& $0.05$	& $0.05$ \\
		Barrier parameter multiplier $\gamma$
			& $2$		& $10$ \\
		Absolute tolerance $\varepsilon_|abs|$
			& $10^{-4}$	& $10^{-7}$ \\
		Relative tolerance $\varepsilon_|rel|$
			& $10^{-4}$	& $10^{-7}$ \\
		\bottomrule
	\end{tabular}
	\vspace{1mm}
	\caption{Algorithm parameters.}
	\label{tab:parameters}
\end{table}

In our experiments we use the signal model \eqref{model}. The frequencies
$\omega_0, \ldots, \omega_{K-1}$ are drawn randomly on $[0,2\pi)$, such that
the minimum separation%
\footnote{The wrap-around distance on $[0,2\pi)$ is used for all frequency differences.}
between any two frequencies is $4\pi/N$.
The coefficients $c_0,\ldots,c_{K-1}$ are generated independently random according to
a circularly symmetric standard complex Gaussian distribution. After generating
the set of $K$ frequencies and coefficients the variance of the noise vector
$\zeta$ is selected such that the desired signal-to-noise ratio (SNR) is
obtained. The regularization parameter $\tau$ is selected from \eqref{tau}
based on the true noise variance.
We assess the algorithms based on their ability to solve AST, which is obtained
by selecting $w=2e_0$ in \eqref{problem}.

We show results for both the L-BFGS and Newton's variants of FastAST%
\footnote{Our code is publicly available at
github.com/thomaslundgaard/fast-ast.}.
For $N\le512$ our implementation uses the Levinson-Durbin
algorithm for Toeplitz inversion, while for $N>512$ it uses the generalized
Schur algorithm where applicable.
The parameters of the algorithm are listed in Table~\ref{tab:parameters}.
It is worth to say a few words about the number of saved difference vectors $M$ in
L-BFGS. On the one hand, selecting larger values of $M$ can decrease the
total number of iterations required (by improving the Hessian approximation), but
on the other hand doing so increases the number of flops required per
iteration. In our numerical experiments we have found that setting it equal to
the size of $u$ ($M=2N-1$) provides a good trade-off. Loosely speaking this
choice allows L-BFGS to perform a full-rank update of the Hessian
approximation, while it does not increase the asymptotic per-iteration
computational complexity.  With this choice the algorithm asymptotically
requires $\calo(N^2)$ flops per iteration.

Performance of the ADMM algorithm%
\footnote{We use the implementation from github.com/badrinarayan/astlinespec.}
\cite{Bhaskar:2013} is also shown along with that of CVX \cite{cvx} applied with both the SeDuMi\cite{Sturm:99} and Mosek\footnote{mosek.com} solvers.

\subsection{Solution Accuracy Per Iteration}

\begin{figure}[t]
		\centering
		\input{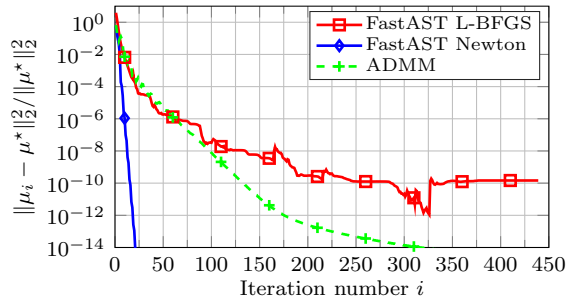}
	\vspace{-3mm}
	\caption{Solution accuracy versus iteration. The signal length is $N=64$,
	the number of sinusoids is $K=6$ and the SNR is $20\,\db$.}
	\label{fig:iters}
\end{figure}

For this investigation a ground-truth solution of \eqref{problem} is obtained
using CVX+SeDuMi with the precision setting set to ``best''. We denote this value as
$\mu^\star$. Fig.~\ref{fig:iters} shows the normalize squared error between $\mu^\star$
and the solution in each iteration of the algorithms. The algorithms ignore the
stopping criteria and run until no further progress can be made towards the
solution.

FastAST Newton converges very fast and a solution of very high accuracy is
obtained within 25 iterations. This is due to the
well-known quadratic convergence of Newton's method. While FastAST
L-BFGS converges significantly slower it requires only $\calo(N^2)$ flops per
iteration versus the $\calo(N^3)$ flops per iteration of FastAST Newton.
We therefore cannot, at this point, conclude which version of FastAST is faster
in practice. Note that ADMM on the other hand requires $\calo(N^3)$ flops per
iteration, the same as FastAST Newton, but requires significantly more
iterations.

It is seen that FastAST L-BFGS seems to not make progress after approx.~300 iterations.
This happens due to numerical challenges in evaluating the L-BFGS search
direction. It is well-known that Woodbury's matrix identity, upon which L-BFGS
is based, has limited numerical stability. For this reason FastAST \mbox{L-BFGS}
is unable to obtain a solution of the same accuracy as the SeDuMi and Mosek solvers.
Despite of this, as seen in the following sections, the solution accuracy of
FastAST L-BFGS is sufficiently high in all cases but those with very high SNR.
The tolerance values of FastAST L-BFGS are selected larger than for FastAST
Newton (Table~\ref{tab:parameters}) because of the mentioned numerical issues
with obtaining a high-accuracy solution.

FastAST Newton does not suffer from this problem and can obtain a solution of
about the same accuracy as SeDuMi and Mosek. ADMM can also obtain a solution of high
accuracy but, as can be seen in Fig.~\ref{fig:iters}, it has slow
convergence starting around iteration number $175$. It therefore takes a large
number of iterations to obtain a solution of the same accuracy as SeDuMi/Mosek or
FastAST Newton.

\pgfplotsset{global_axis_style/.style={
	title style={font=\footnotesize},
	legend columns=1,
	legend style={font=\scriptsize},
	legend style={inner xsep=2pt, inner ysep=1pt, nodes={inner sep=0.8pt}},
	legend style={/tikz/every even column/.append style={column sep=2pt}},
	label style={font=\footnotesize},
	xlabel shift=-3pt,
	ylabel shift=-3pt,
	xticklabel style={font=\footnotesize},
	yticklabel style={font=\footnotesize},
	every axis plot/.append style={line width=1pt}
}} 

\setlength\figureheight{32mm}
\setlength\figurewidth{60mm}
\pgfplotsset{local_axis_style/.style={
	mark phase={10},
	mark repeat={50},
	minor x tick num = 1,
	xminorgrids=true
}}

\pgfplotsset{local_axis_style/.style={
	transpose legend,
	legend columns=3,
	xticklabel={
        \pgfkeys{/pgf/fpu=true}
		\pgfmathparse{exp(\tick)}%
		\pgfmathprintnumber[fixed, precision=0]{\pgfmathresult}
        \pgfkeys{/pgf/fpu=false}
      }
}}
\begin{figure*}[t]
	\begin{minipage}[t]{0.5\linewidth}
		\centering
		\pgfplotsset{local_axis_style/.append style={
			ymin = 7e-4, ymax = 1e-2 }} 
%
\definecolor{mycolor1}{rgb}{1.00000,0.00000,1.00000}%
\definecolor{mycolor2}{rgb}{0.00000,1.00000,1.00000}%
\begin{tikzpicture}

\begin{axis}[%
width=0.951\figurewidth,
height=\figureheight,
at={(0\figurewidth,0\figureheight)},
scale only axis,
xmode=log,
xmin=16,
xmax=2048,
xminorticks=true,
xlabel style={font=\color{white!15!black}},
xlabel={$N$},
ymode=log,
ymin=0.00089584905365534,
ymax=0.00204703994085593,
yminorticks=true,
ylabel style={font=\color{white!15!black}},
ylabel={Reconstruction NMSE},
axis background/.style={fill=white},
xmajorgrids,
xminorgrids,
ymajorgrids,
yminorgrids,
legend style={legend cell align=left, align=left, draw=white!15!black},
global_axis_style, local_axis_style,
xtick = {1.6000e+01,3.2000e+01,6.4000e+01,1.2800e+02,2.5600e+02,5.1200e+02,1.0240e+03,2.0480e+03},
ytick = {1.0e-10,1.0e-09,1.0e-08,1.0e-07,1.0e-06,1.0e-05,1.0e-04,1.0e-03,1.0e-02,1.0e-01,1.0e+00,1.0e+01,1.0e+02,1.0e+03,1.0e+04,1.0e+05},
minor ytick = {1.0e-10,2.0e-10,3.0e-10,4.0e-10,5.0e-10,6.0e-10,7.0e-10,8.0e-10,9.0e-10,1.0e-09,2.0e-09,3.0e-09,4.0e-09,5.0e-09,6.0e-09,7.0e-09,8.0e-09,9.0e-09,1.0e-08,2.0e-08,3.0e-08,4.0e-08,5.0e-08,6.0e-08,7.0e-08,8.0e-08,9.0e-08,1.0e-07,2.0e-07,3.0e-07,4.0e-07,5.0e-07,6.0e-07,7.0e-07,8.0e-07,9.0e-07,1.0e-06,2.0e-06,3.0e-06,4.0e-06,5.0e-06,6.0e-06,7.0e-06,8.0e-06,9.0e-06,1.0e-05,2.0e-05,3.0e-05,4.0e-05,5.0e-05,6.0e-05,7.0e-05,8.0e-05,9.0e-05,1.0e-04,2.0e-04,3.0e-04,4.0e-04,5.0e-04,6.0e-04,7.0e-04,8.0e-04,9.0e-04,1.0e-03,2.0e-03,3.0e-03,4.0e-03,5.0e-03,6.0e-03,7.0e-03,8.0e-03,9.0e-03,1.0e-02,2.0e-02,3.0e-02,4.0e-02,5.0e-02,6.0e-02,7.0e-02,8.0e-02,9.0e-02,1.0e-01,2.0e-01,3.0e-01,4.0e-01,5.0e-01,6.0e-01,7.0e-01,8.0e-01,9.0e-01,1.0e+00,2.0e+00,3.0e+00,4.0e+00,5.0e+00,6.0e+00,7.0e+00,8.0e+00,9.0e+00,1.0e+01,2.0e+01,3.0e+01,4.0e+01,5.0e+01,6.0e+01,7.0e+01,8.0e+01,9.0e+01,1.0e+02,2.0e+02,3.0e+02,4.0e+02,5.0e+02,6.0e+02,7.0e+02,8.0e+02,9.0e+02,1.0e+03,2.0e+03,3.0e+03,4.0e+03,5.0e+03,6.0e+03,7.0e+03,8.0e+03,9.0e+03,1.0e+04,2.0e+04,3.0e+04,4.0e+04,5.0e+04,6.0e+04,7.0e+04,8.0e+04,9.0e+04,1.0e+05,2.0e+05,3.0e+05,4.0e+05,5.0e+05,6.0e+05,7.0e+05,8.0e+05,9.0e+05}
]
\addplot [color=red, mark=square, mark options={solid, red}]
  table[row sep=crcr]{%
16	0.0020369925399371\\
32	0.00146086270670006\\
64	0.00153924658863511\\
128	0.00172306328187387\\
256	0.00171025409965305\\
512	0.00169260371390043\\
1024	0.00171796275319855\\
2048	0.00173986902898587\\
};
\addlegendentry{FastAST L-BFGS}

\addplot [color=blue, mark=diamond, mark options={solid, blue}]
  table[row sep=crcr]{%
16	0.00204675087675347\\
32	0.00145764459154625\\
64	0.00153782643706066\\
128	0.00172211242837243\\
256	0.0017100810404626\\
512	0.00168909097950233\\
1024	0.00171560105210941\\
};
\addlegendentry{FastAST Newton}

\addplot [color=green, dashed, mark=+, mark options={solid, green}]
  table[row sep=crcr]{%
16	0.00204703994085593\\
32	0.0014575103937561\\
64	0.00153801735682853\\
128	0.00172224837800377\\
256	0.00171044330605312\\
512	0.00169669257898609\\
};
\addlegendentry{ADMM}

\addplot [color=mycolor1, dotted, mark=x, mark options={solid, mycolor1}]
  table[row sep=crcr]{%
16	0.00204675087675347\\
32	0.00145763826594231\\
64	0.00153782643706066\\
128	0.00172204855915906\\
256	0.0017100810404626\\
};
\addlegendentry{CVX+SeDuMi}

\addplot [color=mycolor2, dotted, mark=triangle, mark options={solid, rotate=180, mycolor2}]
  table[row sep=crcr]{%
16	0.00204675087675347\\
32	0.00145764459154625\\
64	0.00153782643706066\\
128	0.00172209522068311\\
256	0.0017100810404626\\
};
\addlegendentry{CVX+Mosek}

\addplot [color=black, dashdotted]
  table[row sep=crcr]{%
16	0.00125906451186677\\
32	0.000973195008595281\\
64	0.00089584905365534\\
128	0.00100628100986409\\
256	0.000988784681145047\\
512	0.000991914381941318\\
1024	0.0010097766687906\\
2048	0.00100406890793736\\
};
\addlegendentry{Oracle}

\end{axis}
\end{tikzpicture}%
	\end{minipage}%
	\begin{minipage}[t]{0.5\linewidth}
		\centering
		\pgfplotsset{local_axis_style/.append style={
			ymin = 1e-10, ymax = 1e-5 }} 
%
\definecolor{mycolor1}{rgb}{1.00000,0.00000,1.00000}%
\definecolor{mycolor2}{rgb}{0.00000,1.00000,1.00000}%
\begin{tikzpicture}

\begin{axis}[%
width=0.951\figurewidth,
height=\figureheight,
at={(0\figurewidth,0\figureheight)},
scale only axis,
xmode=log,
xmin=16,
xmax=2048,
xminorticks=true,
xlabel style={font=\color{white!15!black}},
xlabel={$N$},
ymode=log,
ymin=1e-10,
ymax=1e-05,
yminorticks=true,
ylabel style={font=\color{white!15!black}},
ylabel={Conditional Freq. MSE},
axis background/.style={fill=white},
xmajorgrids,
xminorgrids,
ymajorgrids,
yminorgrids,
global_axis_style, local_axis_style,
xtick = {1.6000e+01,3.2000e+01,6.4000e+01,1.2800e+02,2.5600e+02,5.1200e+02,1.0240e+03,2.0480e+03},
ytick = {1.0e-10,1.0e-09,1.0e-08,1.0e-07,1.0e-06,1.0e-05,1.0e-04,1.0e-03,1.0e-02,1.0e-01,1.0e+00,1.0e+01,1.0e+02,1.0e+03,1.0e+04,1.0e+05},
minor ytick = {1.0e-10,2.0e-10,3.0e-10,4.0e-10,5.0e-10,6.0e-10,7.0e-10,8.0e-10,9.0e-10,1.0e-09,2.0e-09,3.0e-09,4.0e-09,5.0e-09,6.0e-09,7.0e-09,8.0e-09,9.0e-09,1.0e-08,2.0e-08,3.0e-08,4.0e-08,5.0e-08,6.0e-08,7.0e-08,8.0e-08,9.0e-08,1.0e-07,2.0e-07,3.0e-07,4.0e-07,5.0e-07,6.0e-07,7.0e-07,8.0e-07,9.0e-07,1.0e-06,2.0e-06,3.0e-06,4.0e-06,5.0e-06,6.0e-06,7.0e-06,8.0e-06,9.0e-06,1.0e-05,2.0e-05,3.0e-05,4.0e-05,5.0e-05,6.0e-05,7.0e-05,8.0e-05,9.0e-05,1.0e-04,2.0e-04,3.0e-04,4.0e-04,5.0e-04,6.0e-04,7.0e-04,8.0e-04,9.0e-04,1.0e-03,2.0e-03,3.0e-03,4.0e-03,5.0e-03,6.0e-03,7.0e-03,8.0e-03,9.0e-03,1.0e-02,2.0e-02,3.0e-02,4.0e-02,5.0e-02,6.0e-02,7.0e-02,8.0e-02,9.0e-02,1.0e-01,2.0e-01,3.0e-01,4.0e-01,5.0e-01,6.0e-01,7.0e-01,8.0e-01,9.0e-01,1.0e+00,2.0e+00,3.0e+00,4.0e+00,5.0e+00,6.0e+00,7.0e+00,8.0e+00,9.0e+00,1.0e+01,2.0e+01,3.0e+01,4.0e+01,5.0e+01,6.0e+01,7.0e+01,8.0e+01,9.0e+01,1.0e+02,2.0e+02,3.0e+02,4.0e+02,5.0e+02,6.0e+02,7.0e+02,8.0e+02,9.0e+02,1.0e+03,2.0e+03,3.0e+03,4.0e+03,5.0e+03,6.0e+03,7.0e+03,8.0e+03,9.0e+03,1.0e+04,2.0e+04,3.0e+04,4.0e+04,5.0e+04,6.0e+04,7.0e+04,8.0e+04,9.0e+04,1.0e+05,2.0e+05,3.0e+05,4.0e+05,5.0e+05,6.0e+05,7.0e+05,8.0e+05,9.0e+05}
]
\addplot [color=red, mark=square, mark options={solid, red}, forget plot]
  table[row sep=crcr]{%
16	2.18207005792844e-06\\
32	5.58720399763362e-07\\
64	1.93901095327913e-07\\
128	3.72991459809485e-08\\
256	1.15320546487081e-08\\
512	3.1927736494928e-09\\
1024	8.62823933976179e-10\\
2048	1.50562766727149e-10\\
};
\addplot [color=blue, mark=diamond, mark options={solid, blue}, forget plot]
  table[row sep=crcr]{%
16	2.19034649403063e-06\\
32	5.56608515837467e-07\\
64	1.93996756738557e-07\\
128	3.72723295982925e-08\\
256	1.15337462244683e-08\\
512	3.19303076057873e-09\\
1024	8.60792566304472e-10\\
};
\addplot [color=green, dashed, mark=+, mark options={solid, green}, forget plot]
  table[row sep=crcr]{%
16	2.19416279148421e-06\\
32	5.56670435557068e-07\\
64	1.94059915067549e-07\\
128	3.72880555304758e-08\\
256	1.15464235463066e-08\\
512	3.18230035462247e-09\\
};
\addplot [color=mycolor1, dotted, mark=x, mark options={solid, mycolor1}, forget plot]
  table[row sep=crcr]{%
16	2.19034649403063e-06\\
32	5.56606784414083e-07\\
64	1.93996756738557e-07\\
128	3.72764452674614e-08\\
256	1.15337462244683e-08\\
};
\addplot [color=mycolor2, dotted, mark=triangle, mark options={solid, rotate=180, mycolor2}, forget plot]
  table[row sep=crcr]{%
16	2.19034649403063e-06\\
32	5.56608515837467e-07\\
64	1.93996756738557e-07\\
128	3.72723295982925e-08\\
256	1.15337462244683e-08\\
};
\end{axis}
\end{tikzpicture}%
	\end{minipage}
	\begin{minipage}[t]{0.5\linewidth}
		\centering
		\pgfplotsset{local_axis_style/.append style={
			ymin =, ymax = }} 
%
\definecolor{mycolor1}{rgb}{1.00000,0.00000,1.00000}%
\definecolor{mycolor2}{rgb}{0.00000,1.00000,1.00000}%
\begin{tikzpicture}

\begin{axis}[%
width=0.951\figurewidth,
height=\figureheight,
at={(0\figurewidth,0\figureheight)},
scale only axis,
unbounded coords=jump,
xmode=log,
xmin=16,
xmax=2048,
xminorticks=true,
xlabel style={font=\color{white!15!black}},
xlabel={$N$},
ymode=log,
ymin=10,
ymax=1000,
yminorticks=true,
ylabel style={font=\color{white!15!black}},
ylabel={Iterations},
axis background/.style={fill=white},
xmajorgrids,
xminorgrids,
ymajorgrids,
yminorgrids,
global_axis_style, local_axis_style,
xtick = {1.6000e+01,3.2000e+01,6.4000e+01,1.2800e+02,2.5600e+02,5.1200e+02,1.0240e+03,2.0480e+03},
ytick = {1.0e-10,1.0e-09,1.0e-08,1.0e-07,1.0e-06,1.0e-05,1.0e-04,1.0e-03,1.0e-02,1.0e-01,1.0e+00,1.0e+01,1.0e+02,1.0e+03,1.0e+04,1.0e+05},
minor ytick = {1.0e-10,2.0e-10,3.0e-10,4.0e-10,5.0e-10,6.0e-10,7.0e-10,8.0e-10,9.0e-10,1.0e-09,2.0e-09,3.0e-09,4.0e-09,5.0e-09,6.0e-09,7.0e-09,8.0e-09,9.0e-09,1.0e-08,2.0e-08,3.0e-08,4.0e-08,5.0e-08,6.0e-08,7.0e-08,8.0e-08,9.0e-08,1.0e-07,2.0e-07,3.0e-07,4.0e-07,5.0e-07,6.0e-07,7.0e-07,8.0e-07,9.0e-07,1.0e-06,2.0e-06,3.0e-06,4.0e-06,5.0e-06,6.0e-06,7.0e-06,8.0e-06,9.0e-06,1.0e-05,2.0e-05,3.0e-05,4.0e-05,5.0e-05,6.0e-05,7.0e-05,8.0e-05,9.0e-05,1.0e-04,2.0e-04,3.0e-04,4.0e-04,5.0e-04,6.0e-04,7.0e-04,8.0e-04,9.0e-04,1.0e-03,2.0e-03,3.0e-03,4.0e-03,5.0e-03,6.0e-03,7.0e-03,8.0e-03,9.0e-03,1.0e-02,2.0e-02,3.0e-02,4.0e-02,5.0e-02,6.0e-02,7.0e-02,8.0e-02,9.0e-02,1.0e-01,2.0e-01,3.0e-01,4.0e-01,5.0e-01,6.0e-01,7.0e-01,8.0e-01,9.0e-01,1.0e+00,2.0e+00,3.0e+00,4.0e+00,5.0e+00,6.0e+00,7.0e+00,8.0e+00,9.0e+00,1.0e+01,2.0e+01,3.0e+01,4.0e+01,5.0e+01,6.0e+01,7.0e+01,8.0e+01,9.0e+01,1.0e+02,2.0e+02,3.0e+02,4.0e+02,5.0e+02,6.0e+02,7.0e+02,8.0e+02,9.0e+02,1.0e+03,2.0e+03,3.0e+03,4.0e+03,5.0e+03,6.0e+03,7.0e+03,8.0e+03,9.0e+03,1.0e+04,2.0e+04,3.0e+04,4.0e+04,5.0e+04,6.0e+04,7.0e+04,8.0e+04,9.0e+04,1.0e+05,2.0e+05,3.0e+05,4.0e+05,5.0e+05,6.0e+05,7.0e+05,8.0e+05,9.0e+05}
]
\addplot [color=red, mark=square, mark options={solid, red}, forget plot]
  table[row sep=crcr]{%
16	68.22\\
32	89.97\\
64	143.77\\
128	204.14\\
256	273.93\\
512	340.37\\
1024	439.86\\
2048	558.26\\
};
\addplot [color=blue, mark=diamond, mark options={solid, blue}, forget plot]
  table[row sep=crcr]{%
16	21.72\\
32	22.09\\
64	22.52\\
128	22.69\\
256	22.93\\
512	23.08\\
1024	23.57\\
};
\addplot [color=green, dashed, mark=+, mark options={solid, green}, forget plot]
  table[row sep=crcr]{%
16	27.64\\
32	45.87\\
64	77.53\\
128	144.09\\
256	277.47\\
512	523.98\\
};
\addplot [color=mycolor1, dotted, mark=x, mark options={solid, mycolor1}, forget plot]
  table[row sep=crcr]{%
16	nan\\
32	nan\\
64	nan\\
128	nan\\
256	nan\\
};
\addplot [color=mycolor2, dotted, mark=triangle, mark options={solid, rotate=180, mycolor2}, forget plot]
  table[row sep=crcr]{%
16	nan\\
32	nan\\
64	nan\\
128	nan\\
256	nan\\
};
\end{axis}
\end{tikzpicture}%
	\end{minipage}%
	\begin{minipage}[t]{0.5\linewidth}
		\centering
		\pgfplotsset{local_axis_style/.append style={
			ymin = 4e-3, ymax = 1e2 }} 
		\input{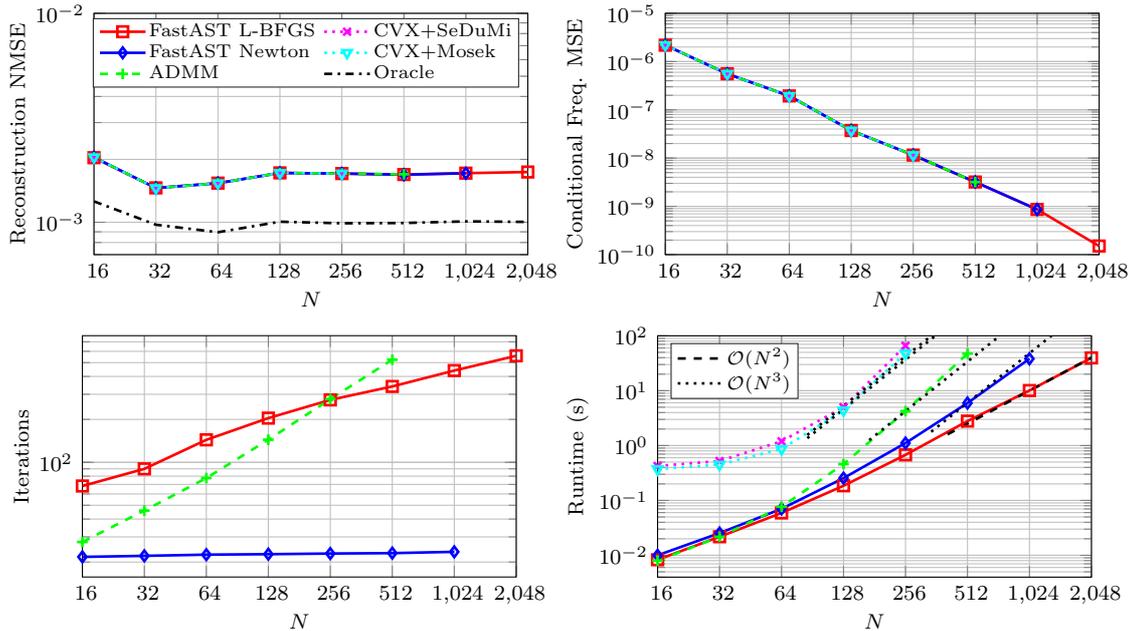}
	\end{minipage}%
	\vspace{-3mm}
	\caption{Simulation results for varying problem size $N$. The SNR is
	$20\,\db$ and the number of sinusoids $K$ is selected as $N/10$ rounded to
	the nearest integer. Results are averaged over $100$ Monte Carlo trials.
	The legend applies to all plots; only the NMSE of Oracle is shown. In the
	figure with runtime the asymptotic per-iteration computational complexity is also
	plotted.}
	\label{fig:N}
\end{figure*}

\subsection{Metrics}
In the following we perform a Monte Carlo simulation study.  Four metrics of
algorithm performance and behaviour are considered: normalized mean-square
error (NMSE) of the reconstructed signal $x$; mean-square error (MSE) of the
frequencies $\{\omega_k\}$ conditioned on successful recovery; number of
iterations and algorithm runtime.  The NMSE of the reconstructed signal is
obtained by estimating the frequencies from the dual polynomial as described in
\cite{Bhaskar:2013} and using these to obtain the least-squares solution for
the coefficients. An estimate of $x$ is then obtained by inserting into
\eqref{model}. This estimate is also known as the \textit{debiased} solution
and it is known to have smaller NMSE than the estimate of $x$ directly obtained
as the solution of \eqref{problem} \cite{Bhaskar:2013}.  In the evaluation of
the signal reconstruction the performance of an Oracle estimator is also shown.
The Oracle estimator knows the true frequencies and estimates the coefficients
using least-squares.

To directly assess the accuracy with which the frequencies are estimated we
present the MSE of the frequency estimates obtained from the dual polynomial. 
The MSE of the frequency estimates is only calculated based on these Monto
Carlo trails in which the set of frequencies is successfully recovered.
Successful recovery is understood as
correct estimation of the model order $K$ and that
all frequency estimates are within a distance of $\pi/N$ from their true value.
The association of the estimated to the true frequencies is obtained by
minimizing the frequency MSE using the Hungarian method \cite{kuhn-hungarian}.

The simulations are performed on a T470p Lenovo, with an Intel(R)
Core(TM) i7-7820HQ CPU @ 2.90GHz, using MATLAB R2018b. MATLAB is
restricted to only use a single CPU core, such that the runtime of the
algorithms can be compared without differences in the parallelism
achieved in the implementations. The computationally heavy steps of
FastAST and ADMM are implemented in native code using the automatic
code generation (``codegen'') feature of MATLAB.

\subsection{Performance Versus Problem Size}

The performance versus problem size $N$ is depicted in Fig.~\ref{fig:N}. First
note that all algorithms give the same estimation accuracy at all problem
sizes, providing strong evidence that they correctly solve \eqref{problem}.

The number of iterations of FastAST L-BFGS increases with $N$. It is
then expected that the total runtime asymptotically scales at a rate
above the per-iteration cost of $\calo(N^2)$ flops. Even still, the runtime for
$N$ up to $2,048$ scales at a rate of about $\calo(N^2)$.
The number of iterations of FastAST Newton is practically independent of $N$.
We then
expect the total runtime to scale asymptotically as $\calo(N^3)$. In practice
it scales a little better for the values of $N$ considered here. The number of
iterations of ADMM increases significantly with $N$ (doubling $N$ roughly
doubles the number of iterations). This in turn means that
the runtime scales faster than the asymptotic per-iteration cost of
$\calo(N^3)$ flops.

In conclusion both variants of FastAST are faster than ADMM already at $N=128$
and their runtime scales at a rate much slower than ADMM. This means that they are
significantly faster than ADMM for large values of $N$. For large $N$ it is
also clear that the L-BFGS variant of FastAST is significantly faster than the
Newton variant.

\subsection{Performance Versus Signal-to-Noise Ratio}
Fig.~\ref{fig:snr} shows performance versus the SNR level. Note that the
conditional MSE of the frequency estimates is not shown for $0\,\db$ SNR
because there are no Monte Carlo trials with successful recovery of the
frequencies at this SNR.

At SNR up to $30\,\db$ all the algorithms perform the same in terms of NMSE of
$x$ and conditional MSE of the frequency estimates. This means that all
algorithms have found a sufficiently accurate solution of \eqref{problem}
(relative to the SNR). In SNR larger than $30\,\db$ FastAST L-BFGS shows a
degraded solution accuracy compared to the remaining algorithms. This is
due to the mentioned numerical issues and the consequently larger
tolerances selected (cf. Table~\ref{tab:parameters}).

In terms of number of iterations and runtime note that both variants of
FastAST show roughly unchanged behaviour with different SNR. ADMM on the other hand
requires more iterations and has larger runtime for large SNR.
In large SNR it is evident that FastAST Newton is preferred due to lower
runtime than ADMM and higher estimation accuracy than FastAST L-BFGS.

\pgfplotsset{local_axis_style/.style={
	transpose legend,
	legend columns=3,
}}
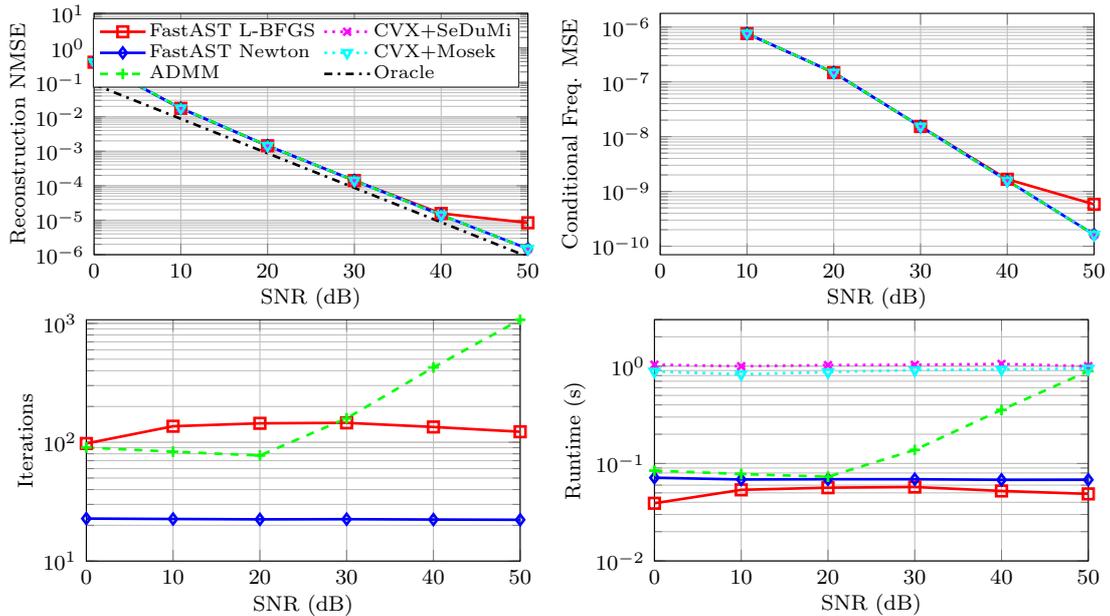
\begin{figure*}[t]
	\begin{minipage}[t]{0.5\linewidth}
		\centering
		\pgfplotsset{local_axis_style/.append style={
			ymin = 1e-6, ymax = 1e1 }} 
%
\definecolor{mycolor1}{rgb}{1.00000,0.00000,1.00000}%
\definecolor{mycolor2}{rgb}{0.00000,1.00000,1.00000}%
\begin{tikzpicture}

\begin{axis}[%
width=0.951\figurewidth,
height=\figureheight,
at={(0\figurewidth,0\figureheight)},
scale only axis,
xmin=0,
xmax=50,
xlabel style={font=\color{white!15!black}},
xlabel={SNR (dB)},
ymode=log,
ymin=1e-07,
ymax=1,
yminorticks=true,
ylabel style={font=\color{white!15!black}},
ylabel={Reconstruction NMSE},
axis background/.style={fill=white},
xmajorgrids,
ymajorgrids,
yminorgrids,
legend style={legend cell align=left, align=left, draw=white!15!black},
global_axis_style, local_axis_style,
xtick = {0.0000e+00,1.0000e+01,2.0000e+01,3.0000e+01,4.0000e+01,5.0000e+01},
ytick = {1.0e-10,1.0e-09,1.0e-08,1.0e-07,1.0e-06,1.0e-05,1.0e-04,1.0e-03,1.0e-02,1.0e-01,1.0e+00,1.0e+01,1.0e+02,1.0e+03,1.0e+04,1.0e+05},
minor ytick = {1.0e-10,2.0e-10,3.0e-10,4.0e-10,5.0e-10,6.0e-10,7.0e-10,8.0e-10,9.0e-10,1.0e-09,2.0e-09,3.0e-09,4.0e-09,5.0e-09,6.0e-09,7.0e-09,8.0e-09,9.0e-09,1.0e-08,2.0e-08,3.0e-08,4.0e-08,5.0e-08,6.0e-08,7.0e-08,8.0e-08,9.0e-08,1.0e-07,2.0e-07,3.0e-07,4.0e-07,5.0e-07,6.0e-07,7.0e-07,8.0e-07,9.0e-07,1.0e-06,2.0e-06,3.0e-06,4.0e-06,5.0e-06,6.0e-06,7.0e-06,8.0e-06,9.0e-06,1.0e-05,2.0e-05,3.0e-05,4.0e-05,5.0e-05,6.0e-05,7.0e-05,8.0e-05,9.0e-05,1.0e-04,2.0e-04,3.0e-04,4.0e-04,5.0e-04,6.0e-04,7.0e-04,8.0e-04,9.0e-04,1.0e-03,2.0e-03,3.0e-03,4.0e-03,5.0e-03,6.0e-03,7.0e-03,8.0e-03,9.0e-03,1.0e-02,2.0e-02,3.0e-02,4.0e-02,5.0e-02,6.0e-02,7.0e-02,8.0e-02,9.0e-02,1.0e-01,2.0e-01,3.0e-01,4.0e-01,5.0e-01,6.0e-01,7.0e-01,8.0e-01,9.0e-01,1.0e+00,2.0e+00,3.0e+00,4.0e+00,5.0e+00,6.0e+00,7.0e+00,8.0e+00,9.0e+00,1.0e+01,2.0e+01,3.0e+01,4.0e+01,5.0e+01,6.0e+01,7.0e+01,8.0e+01,9.0e+01,1.0e+02,2.0e+02,3.0e+02,4.0e+02,5.0e+02,6.0e+02,7.0e+02,8.0e+02,9.0e+02,1.0e+03,2.0e+03,3.0e+03,4.0e+03,5.0e+03,6.0e+03,7.0e+03,8.0e+03,9.0e+03,1.0e+04,2.0e+04,3.0e+04,4.0e+04,5.0e+04,6.0e+04,7.0e+04,8.0e+04,9.0e+04,1.0e+05,2.0e+05,3.0e+05,4.0e+05,5.0e+05,6.0e+05,7.0e+05,8.0e+05,9.0e+05}
]
\addplot [color=red, mark=square, mark options={solid, red}]
  table[row sep=crcr]{%
0	0.382759364296716\\
10	0.0175889972645114\\
20	0.00144157876752318\\
30	0.000141156077363516\\
40	1.56997644019249e-05\\
50	8.4184386186294e-06\\
};
\addlegendentry{FastAST L-BFGS}

\addplot [color=blue, mark=diamond, mark options={solid, blue}]
  table[row sep=crcr]{%
0	0.382773541031254\\
10	0.0175914568949126\\
20	0.00144138780823432\\
30	0.000141232462183049\\
40	1.41552466345976e-05\\
50	1.48170186516114e-06\\
};
\addlegendentry{FastAST Newton}

\addplot [color=green, dashed, mark=+, mark options={solid, green}]
  table[row sep=crcr]{%
0	0.386122151500559\\
10	0.0179696631156253\\
20	0.00144162214162272\\
30	0.000141231862801742\\
40	1.41552466345976e-05\\
50	1.48011673403515e-06\\
};
\addlegendentry{ADMM}

\addplot [color=mycolor1, dotted, mark=x, mark options={solid, mycolor1}]
  table[row sep=crcr]{%
0	0.382773541031254\\
10	0.0175914568949126\\
20	0.00144138780823432\\
30	0.000141232462183049\\
40	1.41552466345976e-05\\
50	1.47651641119832e-06\\
};
\addlegendentry{CVX+SeDuMi}

\addplot [color=mycolor2, dotted, mark=triangle, mark options={solid, rotate=180, mycolor2}]
  table[row sep=crcr]{%
0	0.382773541031254\\
10	0.0175914042106102\\
20	0.00144138780823432\\
30	0.000141230308642125\\
40	1.41552466345976e-05\\
50	1.47950171172085e-06\\
};
\addlegendentry{CVX+Mosek}

\addplot [color=black, dashdotted]
  table[row sep=crcr]{%
0	0.0872325627652628\\
10	0.00872325627652627\\
20	0.000872325627652628\\
30	8.72325627652626e-05\\
40	8.72325627652625e-06\\
50	8.72325627652633e-07\\
};
\addlegendentry{Oracle}

\end{axis}
\end{tikzpicture}%
	\end{minipage}%
	\begin{minipage}[t]{0.5\linewidth}
		\centering
		\pgfplotsset{local_axis_style/.append style={
			ymin =, ymax = }} 
%
\definecolor{mycolor1}{rgb}{1.00000,0.00000,1.00000}%
\definecolor{mycolor2}{rgb}{0.00000,1.00000,1.00000}%
\begin{tikzpicture}

\begin{axis}[%
width=0.951\figurewidth,
height=\figureheight,
at={(0\figurewidth,0\figureheight)},
scale only axis,
unbounded coords=jump,
xmin=0,
xmax=50,
xlabel style={font=\color{white!15!black}},
xlabel={SNR (dB)},
ymode=log,
ymin=1e-10,
ymax=1e-06,
yminorticks=true,
ylabel style={font=\color{white!15!black}},
ylabel={Conditional Freq. MSE},
axis background/.style={fill=white},
xmajorgrids,
ymajorgrids,
yminorgrids,
global_axis_style, local_axis_style,
xtick = {0.0000e+00,1.0000e+01,2.0000e+01,3.0000e+01,4.0000e+01,5.0000e+01},
ytick = {1.0e-10,1.0e-09,1.0e-08,1.0e-07,1.0e-06,1.0e-05,1.0e-04,1.0e-03,1.0e-02,1.0e-01,1.0e+00,1.0e+01,1.0e+02,1.0e+03,1.0e+04,1.0e+05},
minor ytick = {1.0e-10,2.0e-10,3.0e-10,4.0e-10,5.0e-10,6.0e-10,7.0e-10,8.0e-10,9.0e-10,1.0e-09,2.0e-09,3.0e-09,4.0e-09,5.0e-09,6.0e-09,7.0e-09,8.0e-09,9.0e-09,1.0e-08,2.0e-08,3.0e-08,4.0e-08,5.0e-08,6.0e-08,7.0e-08,8.0e-08,9.0e-08,1.0e-07,2.0e-07,3.0e-07,4.0e-07,5.0e-07,6.0e-07,7.0e-07,8.0e-07,9.0e-07,1.0e-06,2.0e-06,3.0e-06,4.0e-06,5.0e-06,6.0e-06,7.0e-06,8.0e-06,9.0e-06,1.0e-05,2.0e-05,3.0e-05,4.0e-05,5.0e-05,6.0e-05,7.0e-05,8.0e-05,9.0e-05,1.0e-04,2.0e-04,3.0e-04,4.0e-04,5.0e-04,6.0e-04,7.0e-04,8.0e-04,9.0e-04,1.0e-03,2.0e-03,3.0e-03,4.0e-03,5.0e-03,6.0e-03,7.0e-03,8.0e-03,9.0e-03,1.0e-02,2.0e-02,3.0e-02,4.0e-02,5.0e-02,6.0e-02,7.0e-02,8.0e-02,9.0e-02,1.0e-01,2.0e-01,3.0e-01,4.0e-01,5.0e-01,6.0e-01,7.0e-01,8.0e-01,9.0e-01,1.0e+00,2.0e+00,3.0e+00,4.0e+00,5.0e+00,6.0e+00,7.0e+00,8.0e+00,9.0e+00,1.0e+01,2.0e+01,3.0e+01,4.0e+01,5.0e+01,6.0e+01,7.0e+01,8.0e+01,9.0e+01,1.0e+02,2.0e+02,3.0e+02,4.0e+02,5.0e+02,6.0e+02,7.0e+02,8.0e+02,9.0e+02,1.0e+03,2.0e+03,3.0e+03,4.0e+03,5.0e+03,6.0e+03,7.0e+03,8.0e+03,9.0e+03,1.0e+04,2.0e+04,3.0e+04,4.0e+04,5.0e+04,6.0e+04,7.0e+04,8.0e+04,9.0e+04,1.0e+05,2.0e+05,3.0e+05,4.0e+05,5.0e+05,6.0e+05,7.0e+05,8.0e+05,9.0e+05}
]
\addplot [color=red, mark=square, mark options={solid, red}, forget plot]
  table[row sep=crcr]{%
0	nan\\
10	7.60100139842468e-07\\
20	1.47071551588487e-07\\
30	1.5321094093039e-08\\
40	1.66553431143325e-09\\
50	5.81256635849113e-10\\
};
\addplot [color=blue, mark=diamond, mark options={solid, blue}, forget plot]
  table[row sep=crcr]{%
0	nan\\
10	7.60211566378876e-07\\
20	1.47125077221806e-07\\
30	1.52813832520969e-08\\
40	1.57611655676439e-09\\
50	1.63170347756168e-10\\
};
\addplot [color=green, dashed, mark=+, mark options={solid, green}, forget plot]
  table[row sep=crcr]{%
0	nan\\
10	7.62179117001865e-07\\
20	1.47187967539055e-07\\
30	1.52824722577696e-08\\
40	1.57611655676439e-09\\
50	1.62949663124267e-10\\
};
\addplot [color=mycolor1, dotted, mark=x, mark options={solid, mycolor1}, forget plot]
  table[row sep=crcr]{%
0	nan\\
10	7.60211566378876e-07\\
20	1.47125077221806e-07\\
30	1.52813832520969e-08\\
40	1.57611655676439e-09\\
50	1.62999591336641e-10\\
};
\addplot [color=mycolor2, dotted, mark=triangle, mark options={solid, rotate=180, mycolor2}, forget plot]
  table[row sep=crcr]{%
0	nan\\
10	7.6019998247157e-07\\
20	1.47125077221806e-07\\
30	1.52810841292936e-08\\
40	1.57611655676439e-09\\
50	1.62932142396342e-10\\
};
\end{axis}
\end{tikzpicture}%
	\end{minipage}\vspace{-2mm}
	\begin{minipage}[t]{0.5\linewidth}
		\centering
%
\definecolor{mycolor1}{rgb}{1.00000,0.00000,1.00000}%
\definecolor{mycolor2}{rgb}{0.00000,1.00000,1.00000}%
\begin{tikzpicture}

\begin{axis}[%
width=0.951\figurewidth,
height=\figureheight,
at={(0\figurewidth,0\figureheight)},
scale only axis,
unbounded coords=jump,
xmin=0,
xmax=50,
xlabel style={font=\color{white!15!black}},
xlabel={SNR (dB)},
ymode=log,
ymin=10,
ymax=1073.35,
yminorticks=true,
ylabel style={font=\color{white!15!black}},
ylabel={Iterations},
axis background/.style={fill=white},
xmajorgrids,
ymajorgrids,
yminorgrids,
global_axis_style, local_axis_style,
xtick = {0.0000e+00,1.0000e+01,2.0000e+01,3.0000e+01,4.0000e+01,5.0000e+01},
ytick = {1.0e-10,1.0e-09,1.0e-08,1.0e-07,1.0e-06,1.0e-05,1.0e-04,1.0e-03,1.0e-02,1.0e-01,1.0e+00,1.0e+01,1.0e+02,1.0e+03,1.0e+04,1.0e+05},
minor ytick = {1.0e-10,2.0e-10,3.0e-10,4.0e-10,5.0e-10,6.0e-10,7.0e-10,8.0e-10,9.0e-10,1.0e-09,2.0e-09,3.0e-09,4.0e-09,5.0e-09,6.0e-09,7.0e-09,8.0e-09,9.0e-09,1.0e-08,2.0e-08,3.0e-08,4.0e-08,5.0e-08,6.0e-08,7.0e-08,8.0e-08,9.0e-08,1.0e-07,2.0e-07,3.0e-07,4.0e-07,5.0e-07,6.0e-07,7.0e-07,8.0e-07,9.0e-07,1.0e-06,2.0e-06,3.0e-06,4.0e-06,5.0e-06,6.0e-06,7.0e-06,8.0e-06,9.0e-06,1.0e-05,2.0e-05,3.0e-05,4.0e-05,5.0e-05,6.0e-05,7.0e-05,8.0e-05,9.0e-05,1.0e-04,2.0e-04,3.0e-04,4.0e-04,5.0e-04,6.0e-04,7.0e-04,8.0e-04,9.0e-04,1.0e-03,2.0e-03,3.0e-03,4.0e-03,5.0e-03,6.0e-03,7.0e-03,8.0e-03,9.0e-03,1.0e-02,2.0e-02,3.0e-02,4.0e-02,5.0e-02,6.0e-02,7.0e-02,8.0e-02,9.0e-02,1.0e-01,2.0e-01,3.0e-01,4.0e-01,5.0e-01,6.0e-01,7.0e-01,8.0e-01,9.0e-01,1.0e+00,2.0e+00,3.0e+00,4.0e+00,5.0e+00,6.0e+00,7.0e+00,8.0e+00,9.0e+00,1.0e+01,2.0e+01,3.0e+01,4.0e+01,5.0e+01,6.0e+01,7.0e+01,8.0e+01,9.0e+01,1.0e+02,2.0e+02,3.0e+02,4.0e+02,5.0e+02,6.0e+02,7.0e+02,8.0e+02,9.0e+02,1.0e+03,2.0e+03,3.0e+03,4.0e+03,5.0e+03,6.0e+03,7.0e+03,8.0e+03,9.0e+03,1.0e+04,2.0e+04,3.0e+04,4.0e+04,5.0e+04,6.0e+04,7.0e+04,8.0e+04,9.0e+04,1.0e+05,2.0e+05,3.0e+05,4.0e+05,5.0e+05,6.0e+05,7.0e+05,8.0e+05,9.0e+05}
]
\addplot [color=red, mark=square, mark options={solid, red}, forget plot]
  table[row sep=crcr]{%
0	97.67\\
10	136.42\\
20	144.34\\
30	145.65\\
40	134.49\\
50	122.58\\
};
\addplot [color=blue, mark=diamond, mark options={solid, blue}, forget plot]
  table[row sep=crcr]{%
0	22.79\\
10	22.61\\
20	22.45\\
30	22.52\\
40	22.37\\
50	22.26\\
};
\addplot [color=green, dashed, mark=+, mark options={solid, green}, forget plot]
  table[row sep=crcr]{%
0	90.37\\
10	83.29\\
20	77.35\\
30	158.94\\
40	426.22\\
50	1073.35\\
};
\addplot [color=mycolor1, dotted, mark=x, mark options={solid, mycolor1}, forget plot]
  table[row sep=crcr]{%
0	nan\\
10	nan\\
20	nan\\
30	nan\\
40	nan\\
50	nan\\
};
\addplot [color=mycolor2, dotted, mark=triangle, mark options={solid, rotate=180, mycolor2}, forget plot]
  table[row sep=crcr]{%
0	nan\\
10	nan\\
20	nan\\
30	nan\\
40	nan\\
50	nan\\
};
\end{axis}
\end{tikzpicture}%
	\end{minipage}%
	\begin{minipage}[t]{0.5\linewidth}
		\centering
		\pgfplotsset{local_axis_style/.append style={
			ymin = 1e-2, ymax = 3e0 }} 
%
\definecolor{mycolor1}{rgb}{1.00000,0.00000,1.00000}%
\definecolor{mycolor2}{rgb}{0.00000,1.00000,1.00000}%
\begin{tikzpicture}

\begin{axis}[%
width=0.951\figurewidth,
height=\figureheight,
at={(0\figurewidth,0\figureheight)},
scale only axis,
xmin=0,
xmax=50,
xlabel style={font=\color{white!15!black}},
xlabel={SNR (dB)},
ymode=log,
ymin=0.0390999999999258,
ymax=1.05240000000013,
yminorticks=true,
ylabel style={font=\color{white!15!black}},
ylabel={Runtime (s)},
axis background/.style={fill=white},
xmajorgrids,
ymajorgrids,
yminorgrids,
global_axis_style, local_axis_style,
xtick = {0.0000e+00,1.0000e+01,2.0000e+01,3.0000e+01,4.0000e+01,5.0000e+01},
ytick = {1.0e-10,1.0e-09,1.0e-08,1.0e-07,1.0e-06,1.0e-05,1.0e-04,1.0e-03,1.0e-02,1.0e-01,1.0e+00,1.0e+01,1.0e+02,1.0e+03,1.0e+04,1.0e+05},
minor ytick = {1.0e-10,2.0e-10,3.0e-10,4.0e-10,5.0e-10,6.0e-10,7.0e-10,8.0e-10,9.0e-10,1.0e-09,2.0e-09,3.0e-09,4.0e-09,5.0e-09,6.0e-09,7.0e-09,8.0e-09,9.0e-09,1.0e-08,2.0e-08,3.0e-08,4.0e-08,5.0e-08,6.0e-08,7.0e-08,8.0e-08,9.0e-08,1.0e-07,2.0e-07,3.0e-07,4.0e-07,5.0e-07,6.0e-07,7.0e-07,8.0e-07,9.0e-07,1.0e-06,2.0e-06,3.0e-06,4.0e-06,5.0e-06,6.0e-06,7.0e-06,8.0e-06,9.0e-06,1.0e-05,2.0e-05,3.0e-05,4.0e-05,5.0e-05,6.0e-05,7.0e-05,8.0e-05,9.0e-05,1.0e-04,2.0e-04,3.0e-04,4.0e-04,5.0e-04,6.0e-04,7.0e-04,8.0e-04,9.0e-04,1.0e-03,2.0e-03,3.0e-03,4.0e-03,5.0e-03,6.0e-03,7.0e-03,8.0e-03,9.0e-03,1.0e-02,2.0e-02,3.0e-02,4.0e-02,5.0e-02,6.0e-02,7.0e-02,8.0e-02,9.0e-02,1.0e-01,2.0e-01,3.0e-01,4.0e-01,5.0e-01,6.0e-01,7.0e-01,8.0e-01,9.0e-01,1.0e+00,2.0e+00,3.0e+00,4.0e+00,5.0e+00,6.0e+00,7.0e+00,8.0e+00,9.0e+00,1.0e+01,2.0e+01,3.0e+01,4.0e+01,5.0e+01,6.0e+01,7.0e+01,8.0e+01,9.0e+01,1.0e+02,2.0e+02,3.0e+02,4.0e+02,5.0e+02,6.0e+02,7.0e+02,8.0e+02,9.0e+02,1.0e+03,2.0e+03,3.0e+03,4.0e+03,5.0e+03,6.0e+03,7.0e+03,8.0e+03,9.0e+03,1.0e+04,2.0e+04,3.0e+04,4.0e+04,5.0e+04,6.0e+04,7.0e+04,8.0e+04,9.0e+04,1.0e+05,2.0e+05,3.0e+05,4.0e+05,5.0e+05,6.0e+05,7.0e+05,8.0e+05,9.0e+05}
]
\addplot [color=red, mark=square, mark options={solid, red}, forget plot]
  table[row sep=crcr]{%
0	0.0390999999999258\\
10	0.0537000000001353\\
20	0.0565000000000509\\
30	0.0573999999999069\\
40	0.052299999999741\\
50	0.048799999999901\\
};
\addplot [color=blue, mark=diamond, mark options={solid, blue}, forget plot]
  table[row sep=crcr]{%
0	0.0714999999999782\\
10	0.0685999999997512\\
20	0.0689999999999418\\
30	0.0689999999999418\\
40	0.0680000000002838\\
50	0.0681000000000495\\
};
\addplot [color=green, dashed, mark=+, mark options={solid, green}, forget plot]
  table[row sep=crcr]{%
0	0.0846000000002095\\
10	0.0780000000001382\\
20	0.0732000000002881\\
30	0.137900000000081\\
40	0.355600000000159\\
50	0.886899999999987\\
};
\addplot [color=mycolor1, dotted, mark=x, mark options={solid, mycolor1}, forget plot]
  table[row sep=crcr]{%
0	1.03429999999989\\
10	0.996600000000035\\
20	1.02400000000001\\
30	1.03040000000015\\
40	1.05240000000013\\
50	0.999099999999999\\
};
\addplot [color=mycolor2, dotted, mark=triangle, mark options={solid, rotate=180, mycolor2}, forget plot]
  table[row sep=crcr]{%
0	0.880299999999988\\
10	0.824299999999785\\
20	0.868899999999885\\
30	0.909899999999725\\
40	0.922799999999916\\
50	0.942200000000084\\
};
\end{axis}
\end{tikzpicture}%
	\end{minipage}%
	\vspace{-3mm}
	\caption{Simulation results for varying SNR. The signal length is $N=64$
	and the number of sinusoids is $K=6$. Results are averaged over $100$
	Monte Carlo trials. The legend applies to all plots; only the NMSE of
	Oracle is shown.}
	\label{fig:snr}
\end{figure*}

\section{Conclusions}
The FastAST algorithm presented in this paper provides a fast approach to
solving the atomic norm soft thresholding problem \eqref{problem}. The L-BFGS
variant provides a reasonably accurate solution and is much faster than any
other algorithm for large problem size $N$. If a solution of high accuracy is
requested, which may be desirable in very high SNR, a variant of FastAST based
on Newton's method is also provided. This variant can find a solution of high
accuracy in a small number of iterations. While it is slower than FastAST
L-BFGS, it is significantly faster than the state-of-the-art method based on
ADMM.

The FastAST algorithm is obtained by reformulating the semidefinite program
\eqref{problem} as a non-symmetric conic program
\eqref{conicproblem}. This reformulation is of key importance in obtaining a
fast algorithm. This work has provided an example of an optimization
problem where it is beneficial to formulate it as a non-symmetric conic program
instead of the standard, and much better understood, formulation as a
symmetric conic program.
We have also provided an implementation of a non-symmetric conic solver, thereby
demonstrating the practical feasibility of this class of methods.

We have demonstrated how the L-BFGS two-loop recursion can be modified to allow
a quasi-Newton solution of the barrier problem \eqref{min_g} even when the
barrier parameter $t$ is updated in every iteration. This approach can directly
be applied in other algorithms based on the barrier method, including
primal-only methods.

Finally note that there are many examples of optimization problems of practical
interest which involve a constraint in either the cone of finite
autocorrelation sequences $\calc^*$ or the cone $\calk$. An example is the
gridless SPICE method \cite{Yang:15} for line spectral estimation; or
frequency-domain system identification and filter design as summarized in
\cite{alkire-autocorrelation}. We expect that equally fast primal-dual IPMs can
be derived for all of these problems using the techniques of this paper.
We also expect that it is fairly straight-forward to extend FastAST to atomic
norm minimization with partial observations \cite{Tang:2013} or multiple
measurement vectors \cite{li-multiple}.
An interesting, but less obvious, extension is to the multi-dimensional
harmonic retrieval problem \cite{chi-twodim}; for that purpose the work
\cite{Yang:2016a} may contain some useful insights.

\section*{Acknowledgements}
We would like to thank Lieven Vandenberghe and Martin Skovgaard Andersen for
providing valuable input to the work and pointing us to some important
references.

\appendix
\section{Characterization of \texorpdfstring{$\calk^*$}{K*}}

To characterize the dual cone $\calk^*$ a number of lemmas are needed.


\begin{lemma}
	\label{inequalities_lemma}
	Let $\calk$ be a proper cone and assume $\lambda\ne0$. If $\left<\lambda,\mu\right>\ge0$ for
	every $\mu\in\intr\calk$, then
	$\left<\lambda,\mu\right>\ge0$ for every $\mu\in\calk$.
\end{lemma}
\begin{IEEEproof}
	Let $\tilde\mu\in\calk$ and let $\{\mu_i\}$ be
	a sequence which converges to $\tilde\mu$ with $\mu_i\in\intr\calk$. Then
	$\left<\lambda,\mu_i\right>\ge0$ and so
	$
		\left<\lambda,\tilde\mu\right>
		= \lim_{i\rightarrow\infty} \left<\lambda,\mu_i\right>
		\ge 0,
	$
	completing the proof.
\end{IEEEproof}

\begin{lemma}
	\label{interior_lemma}
	Let $\calk$ be a proper cone. The interior of its dual is given by
	\begin{align}
		\intr\calk^* = \{ \lambda :
		\left<\lambda,\mu\right> > 0 \;\; \forall \; \mu\in\calk\}.
	\end{align}
\end{lemma}
\begin{IEEEproof}
	See \cite{BoVa:04}, exercise 2.31.
\end{IEEEproof}

To formulate the next lemma, the dual barrier of $F$ is introduced:
\begin{align}
	F^*(\lambda) = \sup\,\{-\left<\lambda,\mu\right> - F(\mu) : \mu\in\intr\calk \}.
	\label{dualbarrier}
\end{align}
This function is a slight modification
($-\left<\lambda,\mu\right>$ replaces $\left<\lambda,\mu\right>$)
of the convex conjugate of $F$. It turns out that $F^*$ is a LH barrier
function for the dual cone $\calk^*$ \cite{BoVa:04,NeTo-selfscaled}. Its usefulness for our
purposes lies in the following property.
\begin{lemma}
	\label{dualbarrier_lemma}
	Assume $\lambda\ne0$ and let $\calk$ be a proper cone with corresponding LH
	barrier function $F$. Then $\lambda\in\intr\calk^*$ if and only if
	$F^*(\lambda)<\infty$
	(i.e., $F^*$ is bounded above.)
\end{lemma}
\begin{IEEEproof}
	We first prove the direct implication. Reasoning by contradiction, assume
	that $F^*(\lambda)<\infty$ and that there exists a $\mu\in\intr\calk$ such that
	$\left<\lambda,\mu\right><0$. Then $\alpha\mu\in\intr\calk$ for all $\alpha>0$.
	But
	$\lim_{\alpha\rightarrow\infty} - \left<\lambda,\alpha\mu\right> - F(\alpha\mu) =
	\lim_{\alpha\rightarrow\infty} - \alpha\left<\lambda,\mu\right> - F(\mu) +
	\theta_F\log(\alpha) = \infty$, a contradiction, so
	$\left<\lambda,\mu\right>\ge0$ for every $\mu\in\intr\calk$. By Lemma
	\ref{inequalities_lemma} we have $\left<\lambda,\mu\right>\ge0$ for every
	$\mu\in\calk$, thus $\lambda\in\calk^*$. Since $F^*$ is a LH barrier function for
	$\calk^*$, it is easy to show that $F^*(\lambda)<\infty$ implies
	$\lambda\notin\bdy\calk^*$, so $\lambda\in\intr\calk^*$.

	To prove the converse assume $\lambda\in\intr\calk^*$. Then by Lemma
	\ref{interior_lemma}, we have $\left<\lambda,\mu\right>>0$ for all
	$\mu\in\calk$. It follows that there exists an $\varepsilon>0$ such that
	$\left<\lambda,\tilde\mu\right>\ge\varepsilon$ for every $\tilde\mu\in\calk$ with
	$\norm{\tilde\mu}_2=1$. By continuity of $F$ it can also be shown that there exists a
	$\delta$ such that $F(\tilde\mu)\ge\delta$ for every $\tilde\mu\in\calk$ with
	$\norm{\tilde\mu}_2=1$. With $\tilde\mu=\mu/\norm{\mu}_2$, the objective in \eqref{dualbarrier} obeys
	\begin{align*}
		- \left<\lambda,\mu\right> - F(\mu)
			&= - \norm{\mu}_2 \left<\lambda,\tilde\mu\right>
				- F\!\left(\norm{\mu}_2\tilde\mu\right) \\
			&= - \norm{\mu}_2 \left<\lambda,\tilde\mu\right>
				- F\!\left(\tilde\mu\right) + \theta_F\log(\norm{\mu}_2) \\
			&\le - \norm{\mu}_2 \varepsilon - \delta + \theta_F \log(\norm{\mu}_2).
	\end{align*}
	The second equality follows from logarithmic homogeneity of $F$.
	This function is bounded above and so $F^*(\lambda)<\infty$.
\end{IEEEproof}

We are now ready to give the desired proof.
\begin{IEEEproof}[Proof of Lemma \ref{Kstar_lemma}]
	It is easy to show the following:
	\begin{enumerate}
		\item If $\rho<0$, then $\lambda\notin\calk^*$.
		\item If $\rho=0$ and $s\neq0$, then $\lambda\notin\calk^*$.
		\item If $\rho=0$ and $s=0$, then $\lambda\in\calk^*$ if and only if
			$z\in\calc^*$.
	\end{enumerate}
	The first and second property are shown by constructing a $\mu\in\calk$
	such that $\left<\lambda,\mu\right><0$.
	The third property is shown by writing $\left<\lambda,\mu\right>=z\T u\ge0$
	for all $\mu\in\calk$ if and only if $z\T u\ge0$ for all $u\in\calc$.

	The only case we have not considered so far is $\rho>0$.
	By Lemma \ref{dualbarrier_lemma} and \eqref{F}, we have $\lambda\in\intr\calk^*$ if and only if
	$F^*(\lambda)<\infty$, i.e., when
	\begin{align*}
		h(\mu) = - \rho v - \re(s\h x) - z\T u + \log|T(u)|
			+ \log(v - x\h T\ii(u) x)
	\end{align*}
	is bounded above on the domain $\mu\in\intr\calk$.
	The function $h$ is concave and by setting the gradient equal to zero we
	get optimal points
	\begin{align*}
		v^\star &= \rho\ii + (2\rho)^{-2} s\h T(u) s \\
		x^\star &= - (2\rho)\ii T(u)s.
	\end{align*}
	It is easy to show that if	$u\in\intr\calc$, then
	$(v^\star, x^\star, u)\T\in\intr\calk$.
	Inserting into $h(\mu)$ we obtain
	\begin{align*}
		h(\mu) &\le - z\T u + \frac{1}{4\rho} s\h T(u) s
		 - 1 - \log(\rho) + \log|T(u)| \\
		 &= - c(\lambda)\T u
			- 1 - \log(\rho) + \log|T(u)|,
	\end{align*}
	with $c(\lambda) = z - \frac{1}{4\rho} T^*(ss\h)$.
	For each $u\in\intr\calc$ there exists some corresponding $\mu\in\intr\calk$ such that the above
	holds with equality.

	If $c(\lambda)=0$, the function $h(\mu)$
	is unbounded above on the domain $\mu\in\intr\calk$ and so
	$\lambda\notin\intr\calk^*$.
	If $c(\lambda)\ne0$ we can use Lemma
	\ref{dualbarrier_lemma} because $-\log|T(u)|$ is a LH barrier function for
	$\calc$. So $h(\mu)$ is bounded above on the domain $\mu\in\intr\calk$ if
	and only if $c(\lambda)\in\intr\calc^*$. Tracing back our steps above we have (for
	$\rho>0$) that $\lambda\in\intr\calk^*$ if and only if $c(\lambda)\in\intr\calc^*$.
	Since both of the dual cones are closed sets and $c(\cdot)$ is a continuous
	function, we have $\lambda\in\calk^*$ if and only if $c\in\calc^*$. That
	completes the proof.
\end{IEEEproof}

\section*{References}
\bibliography{latex/IEEEabrv,latex/refs}

\end{document}